\documentclass[12pt]{article}
\newcommand{\kmcomment}[1]{\null} 
\usepackage[a4paper,truedimen,scale={.85,.85},centering]{geometry} 

\usepackage{verbatim,supertabular,theorem,setspace} 
\usepackage{latexsym,amssymb,amsmath,alltt}

\usepackage{mathptmx} 
\usepackage[dvips]{graphicx}
\usepackage[usenames]{color}
\usepackage{url}

{ 
\theorembodyfont{\rmfamily}
\theoremstyle{plain}
\newtheorem{thm}{Theorem}

\newtheorem{defn}{Definition}
\newtheorem{kmProp}{Proposition}
 
\newtheorem{myRemark}{\textbf{Remark}}
}

\newcommand{\ds}{\displaystyle }
\newcommand{\frakS}[1]{\mathfrak{S}_{#1}}
\newcommand{\fraksp}{\mathfrak{s}\mathfrak{p}}
\newcommand{\frakham}{\mathfrak{h}\mathfrak{a}\mathfrak{m}}
\newcommand{\Pkt}[2]{\{#1,#2\}}%
\newcommand{\mR}{\ensuremath{\mathbb{R}}} 
\newcommand{\mZ}{\ensuremath{\mathbb{Z}}} 
\newcommand{\pdel}{\partial} 
\newcommand{\hk}{\hat{k}}
\newcommand{\myd}{d} 

\newcommand{\Z}[2]{z^{#2}_{#1}}
 
\newcommand{\la}{\langle}
\newcommand{\ra}{\rangle}

\newcommand{\CGF}[2]{\text{C}^{#1}_{GF}({\frakham}_4^0,{Sp}(4,\mR))_{#2}} 

\newcommand{\V}[1]{\ensuremath{ V_{\langle #1\rangle}}} 
\newcommand{\RGF}[4]{{\color{black}\overline{\text{C}}}^{#1}_{GF}({\frakham}_{#2}^{#3})_{#4}} 
\newcommand{\CGFF}[4]{\text{C}^{#1}_{GF}({\frakham}_{#2}^{#3})_{#4}} 
\newcommand{\HGFF}[4]{\text{H}^{#1}_{GF}({\frakham}_{#2}^{#3})_{#4}}
\newcommand{\CGFn}[2]{\text{C}^{#1}_{GF}({\frakham}_{2n}^0,{Sp}(2n,\mR))_{#2}}
\newcommand{\myLR}{\ds  \mathop{\cong}^{\text{LR}}} 

\title{Lower weight Gel'fand-Kalinin-Fuks cohomology groups of the formal 
Hamiltonian vector fields on $\mR^4$}

\author{Kentaro Mikami\thanks{
  Akita University, partially supported by Grant-in-Aid for 
  Scientific Research (C) of JSPS, 
 No.26400063, No.23540067 and No.20540059}
\and   
Yasuharu Nakae\thanks{Department of Computer Science and Engineering 
  Akita University} 
  } 
\date{\today} 
\begin{document}
{\allowdisplaybreaks 
\onehalfspacing 
\maketitle 
\renewcommand{\thefootnote}{\fnsymbol{footnote}}
\begin{abstract} 
In this paper, we investigate the relative Gel'fand-Kalinin-Fuks
cohomology groups of the formal Hamiltonian vector fields on
$\ds\mR^{4}$.  
In the case of formal Hamiltonian vector fields on $\ds\mR^{2}$, we computed the relative Gel'fand-Kalinin-Fuks cohomology
groups of weight $<20$ in the paper by Mikami-Nakae-Kodama.  The main
strategy there was decomposing the Gel'fand-Fucks cochain
complex into irreducible factors and picking up the trivial
representations and their concrete bases, and ours is essentially the
same.  
By computer calculation, we determine the relative Gel'fand-Kalinin-Fuks
cohomology groups of the formal Hamiltonian vector fields on
$\ds  \mR^4$ of  weights $2$, $4$ and $6$.  In the case of
weight $2$, the Betti number of the cohomology group is equal to $1$ at
degree $2$ and is $0$ at any other degree.  In weight $4$, the Betti
number is $2$ at degree $4$ and is $0$ at any other degree, and in
weight $6$, the Betti number is $0$ at any degree.  

\noindent 
\textbf{Main revised points:}(April 2014) 

\begin{enumerate}
\item Fix notations of subalgebras of formal Hamiltonian vector fields:
                
\hfil\begin{tabular}{*{5}{|c}|}\hline
older & $\ds\frakham_{2n}^{0}$ & 
$\ds\frakham_{2n}^{1}$ & 
none & 
$\ds\RGF{\bullet}{2n}{1}{w}$ \\ 
now & $\ds\frakham_{2n}^{ }$ & 
$\ds\frakham_{2n}^{0}$ & 
$\ds\frakham_{2n}^{1}$ & 
$\ds\CGFF{\bullet}{2n}{1}{w}$\\\hline\end{tabular} 
\hfil\hfil       

\item $\ds\CGFF{3}{4}{1}{6}$ in p.8 is more simplified. 
               
\item Added two terms $\ds 2\V{7,7}$ to $\ds\CGFF{6}{4}{1}{6}$ in p.11
        (forgotten writing them in the first version).  
\end{enumerate}

\end{abstract}
\section{Introduction} Inspired by \cite{KOT:MORITA}, we are
interested in getting information about the relative Gel'fand-Kalinin-Fuks cohomology groups of
the formal Hamiltonian vector fields on 
$\ds  \mR^{2n}$ of a given 
weight.  

In \cite{M:N:K}, we dealt with the case where $n=1$ of weight $\leq 20$.
In this paper, we investigate the 
relative Gel'fand-Kalinin-Fuks cohomology groups of
the formal Hamiltonian vector fields on 
$\ds  \mR^{4}$.    
Even for $n=1$ or 2, 
the limitation  comes from
overloading of heavy computations of picking up the trivial representations and
their concrete bases.  
Comparing the case where $n=2$ with the case $n=1$,  we encountered more
difficulty of decomposing into irreducible factors of tensor product,
even though the Littlewood-Richardson formula is    theoretically rather
simple. 
So far, the information we have gotten about the relative
Gel'fand-Kalinin-Fuks cohomology groups of the formal Hamiltonian vector
fields on $\mR^4$ is only in the cases of weight=2,4,6,
and the corresponding Euler characteristic numbers are 1,2 and 0,
respectively.  

\section{Splitting cochains by weight} 
We are interested in the standard
linear symplectic space $\ds  \mR^{2n}$. 
The function space $\ds  C^{\infty}(\mR^{2n})$ forms a Lie
algebra with respect to the Poisson bracket $\ds  \Pkt{\cdot}{\cdot}$.   
For the Darboux coordinate $\ds 
(x_1,\ldots,x_n,y_1,\ldots,y_n)$, 
 we have $\ds  \Pkt{x_i}{y_j}= - \Pkt{y_j}{x_i}= \delta_{ij}$
and 
$\ds  \Pkt{x_i}{x_j}=  \Pkt{y_i}{y_j}= 0$.
Then the space of polynomials of $\ds  
x_1,\ldots, y_n$ is a subalgebra with respect to the Poisson bracket.  
The space $\ds\frakham_{2n}^{ }$ 
of Hamiltonian vector fields which have polynomials as Hamiltonian
potentials is a Lie algebra and the map 
$\ds  f \mapsto - H_{f}$ 
from the Hamiltonian potentials to Hamiltonian vector fields 
is a Lie algebra homomorphism with the
kernel $\cong \mR$.  

We look at the
Lie subalgebra $\ds  \frakham_{2n}^{0}$
of $\ds  \frakham_{2n}^{ }$ formed by elements
which vanish at $0$.
$\ds  \frakham_{2n}^{0}$ corresponds to the algebra of polynomials
without linear terms.
 
In this paper, we are interested in the Gel'fand-Kalinin-Fuks cohomology
groups of  $\ds  \frakham_{2n}^{0}$ when $n=2$.   

We can split the polynomial functions by their homogeneity.  The
cochain complex is the exterior algebra of dual of polynomial functions and we
introduce the ``weight`` on the cochain complex as follows:
 
\begin{defn} 
Let $\ds  \frakS{\ell}$ be the dual space of
$\ell$-homogeneous polynomial functions,
and define the weight of each non zero element of
$\ds  \frakS{\ell}$ to be $\ell-2$.  
For each non-zero element of $\ds  
\frakS{\ell_1} \wedge 
\frakS{\ell_2} \wedge \cdots \wedge 
\frakS{\ell_s}$ 
($ \ell_1 \le \ell_2 \le \cdots \le \ell_s$), 
define  its weight to be $\ds  \sum_{i=1}^{s} (\ell_i -2)$. 
\end{defn}

\begin{kmProp}[cf.\cite{KOT:MORITA},\cite{M:N:K}]
The coboundary operator $\myd$ of the Gel'fand-Kalinin-Fuks cochain
complex preserves the weight,
namely, 
if a cochain $\sigma$ is of weight $w$, then 
$\ds  \myd \sigma$ is also of weight $w$.

Hence we can decompose the total space of cochain complex by degree and
weight: namely, 
$$
\CGFF{m}{2n}{ }{w}
 = \text{LinearSpan of }\{  
\sigma \in 
\Lambda^{k_1} 
\frakS{1} \wedge 
\Lambda^{k_2} 
\frakS{2} \wedge \cdots
\wedge 
\Lambda^{k_s} 
\frakS{s} \mid 
\sum_{i=1}^{s} k_i = m \ , \  
\sum_{i=1}^{s} k_i (i-2) = w\ , \ s=1,2,\ldots  
\}$$ 
and we can define the cohomology group 
$ \ds  
\HGFF{m}{2n}{ }{w}
$.  
\end{kmProp}

$\ds  
\CGFF{\bullet}{2n}{0}{w}
$ is the subspace of    
$\ds  
\CGFF{\bullet}{2n}{ }{w}
$ characterized by 
$k_1=0$.
If we restrict our attention to the cochain complex relative to
$Sp(2n,\mR)$,
then it turns out $k_2=0$ (cf.\cite{M:N:K}). 
Thus we first look at the subcomplex
 $\ds  \CGFF{\bullet}{2n}{1}{w}$
of 
$\ds  
\CGFF{\bullet}{2n}{0}{w}
$
spanned by cochains such that $k_1=k_2=0$, where $\ds
\frakham_{2n}^{1}$  is the subalgebra 
the completion of 
$\ds\sum_{i\ge 3} \frakS{i}^{*}$  
defined in \cite{Kont:RW:MR1671725}. 

We consider finite sequences of non-negative integers $(k_3,k_4,\ldots,k_s)$
satisfying 
\begin{equation}
\sum_{i=3}^{s} k_i = m \quad\text{and}\quad 
\sum_{i=3}^{s} k_i (i-2) = w  \label{A:eqn}  \ .  
\end{equation}
Shifting the indices by 2  as  
$\ds  \hk_i = k_{i+2}$ ($i>0$), we see   
\begin{align*}
w &= \hk_1+ 2 \hk_2+ \cdots + t \hk_t \\
&= 
   (\underbrace{t+\cdots+t}_{\hk_t}) + \cdots  
  + (\underbrace{2+\cdots+2}_{\hk_2}) + 
  (\underbrace{1+\cdots+1}_{\hk_1})
= \ell_1+ \ell_2+ \cdots + \ell_m 
\end{align*}
where $t=s-2$ and  $\ell_1 \ge \ell_2 \ge \cdots \ge \ell_{m} \ge 1$.  
This is a partition of $w$ of length  $m$ or a Young
diagram of height $m$  
with $w$ cells.   
Conversely, for a partition of $w$ 
\begin{equation}\label{eqn:B}
\begin{split}
& w = \ell_1+\ell_2+\cdots+\ell_m \\
& \ell_1 \geqq \ell_2 \geqq \cdots \geqq \ell_m \geqq  1
\end{split}
\end{equation}
$ \hk_i = \# \{ j \,|\, \ell_j=i \}$ gives   
a solution  of (\ref{A:eqn}).   
That means there is a one-to-one correspondence between the  
solutions of (\ref{A:eqn}) and all the partitions of $w$ of length $m$ or Young
diagrams of height $m$  
with $w$ cells.

\begin{myRemark}
Since $- \text{Identity}\in Sp(2n,\mR)$, we see that the relative cochain
complex of odd weight must be the zero space, and hence we only deal
with the complexes of even weights.  
\end{myRemark}
\begin{kmProp}\label{exam::a}
When weight=2, 4 or 6,  
the non-trivial cochain groups are as follows: 
\begin{align*} 
\CGFF{1}{2n}{1}{2} & = \frakS{4} ,  
\qquad 
\CGFF{2}{2n}{1}{2} 
 =  \Lambda^2 \frakS{3}  \\
 \\
\CGFF{1}{2n}{1}{4}  & = \frakS{6},  
\qquad 
\CGFF{2}{2n}{1}{4} 
 =  
\left(  \frakS{3} \wedge \frakS{5}  \right)
\oplus \Lambda^2 \frakS{4}  
\cong 
\left(  \frakS{3} \otimes \frakS{5}  \right)
\oplus \Lambda^2 \frakS{4},  
\\   
\CGFF{3}{2n}{1}{4} 
&
 =  \Lambda^{2} \frakS{3} \wedge \frakS{4} 
 \cong  \Lambda^{2} \frakS{3} \otimes    \frakS{4}, \qquad
\CGFF{4}{2n}{1}{4} =  \Lambda^4 \frakS{3}
\intertext{In the above, we identify the exterior product $\frakS{3} \wedge \frakS{5}$
with the tensor product $\frakS{3} \otimes  \frakS{5}$ as vector spaces, 
and we often use this identification without comments.
}
\CGFF{1}{2n}{1}{6} & = \frakS{8} ,  
\qquad 
\CGFF{2}{2n}{1}{6} 
 =  
\left(  \frakS{3} \otimes \frakS{7}  \right)
\oplus \left(  \frakS{4} \otimes \frakS{6}  \right)
\oplus \Lambda^2 \frakS{5}  
\\ 
\CGFF{3}{2n}{1}{6} 
& = \left( \Lambda^2 \frakS{3} \otimes \frakS{6}  \right)
\oplus \left(  \frakS{3} \otimes \frakS{4} \otimes \frakS{5}  \right)
\oplus  \Lambda^{3} \frakS{4}  
\\ 
\CGFF{4}{2n}{1}{6} 
& = \left( \Lambda^3 \frakS{3} \otimes \frakS{5}  \right)
\oplus \left( \Lambda^{2}  \frakS{3} \otimes \Lambda^{2} \frakS{4} \right)
\\
\CGFF{5}{2n}{1}{6} 
& =  \Lambda^4 \frakS{3} \otimes \frakS{4} ,  
\qquad 
\CGFF{6}{2n}{1}{6} 
 = \Lambda^{6} \frakS{3} 
\end{align*}  
In general,  
$\ds  \Lambda^p \frakS{q} = \{\mathbf{0}\}$ if $ p > 
\dim \frakS{q} = (q+2n-1)!/( q! (2n-1)!)$.  
If $n=1$, $\dim \frakS{3} = 4$ and we have 
$
\CGFF{6}{2}{1}{6}=\{\mathbf{0}\}
$.
\end{kmProp}

\textbf{Proof:}
(1) weight$=2$ case: 
When $m=1$, then $\ell_1 = 2$, so we have $\hk_2= 1$ and $\hk_j=0$ ($j
\ne 2$). Thus, $k_4=1$.   
When $m=2$, then $ 2 = \ell_1 + \ell_2$ ( $\ell_1 \ge \ell_2 \ge 1$),  
$\ell_1 = \ell_2 =1$, so we have $\hk_1 = 2$ and 
$\hk_j=0$ ($j \ne 1$). Thus, $k_3 = 2$.  

(2) weight$=4$ case: 
When $m=2$, i.e., $ 4 = \ell_1 + \ell_2$ ( $\ell_1 \ge \ell_2 \ge 1$),
then   
$(\ell_1,\ell_2) = (3,1) $ or  $(2,2)$, 
so we have $(\hk_1 = 1, \hk_3 = 1) $, or   
$\hk_2=2$ . Thus, $(k_3 = 1, k_5=1)$ or $(k_4=2)$. 
When $m=3$, i.e., $ 4 = \ell_1 + \ell_2 +\ell_3$ ( $\ell_1 \ge \ell_2
\ge \ell_3 \ge  1$), $\ell_1=2, \ell_2 = 1, \ell_3=1$. Thus 
$(\hk_2=1, \hk_1=2)$, so $(k_3=2, k_4=1)$.  
When $m=4$, i.e.,
$4=\ell_1+\ell_2+\ell_3+\ell_4$
($\ell_1 \ge \ell_2 \ge \ell_3 \ge \ell_4 \ge 1$),
$\ell_1=\ell_2=\ell_3=\ell_4=1$.
Thus $\hk_1=4$, so $k_3=4$.

(3) weight$=6$ case: Way is the same, we omit the discussion.  

\section{Coboundary operator}
The relative cochain complex is defined by 
$$\ds  
\CGFn{m}{w}
= 
\{ \sigma \in 
\CGFF{m}{2n}{0}{w}
\mid
 i_{\hat{J}(\xi)} \sigma = 0 ,
 i_{\hat{J}(\xi)}  \myd \sigma =
 0 \quad (\forall \xi\in \fraksp(2n,\mR)) \} \ , $$ 
where $\ds  J$ is the momentum mapping of
$Sp(2n,\mR)$ on $\ds  \mR^{2n}$.  
The first condition $ 
 i_{\hat{J}(\xi)} \sigma = 0$ 
 ($\forall \xi\in \fraksp(2n,\mR)$) 
 means that $\ds  \frakS{2}$ does
 never appear.  The second condition 
$ i_{\hat{J}(\xi)}  \myd \sigma =
 0$ ($\forall \xi\in \fraksp(2n,\mR)$) 
means the cochain complex consists of the trivial representation spaces.    

In the rest of this paper, we use the global variables 
$\ds  x_1, x_2, x_3 , x_4$ on $\ds  \mR^4$ 
and the symplectic form 
$\omega$ on $\ds  \mR^4$ is given by 
$\ds 
\omega( \frac{\pdel}{\pdel x_1},\frac{\pdel}{\pdel x_4}) = - 
\omega( \frac{\pdel}{\pdel x_4},\frac{\pdel}{\pdel x_1}) = 1$,  
$\ds 
\omega( \frac{\pdel}{\pdel x_2},\frac{\pdel}{\pdel x_3}) = - 
\omega( \frac{\pdel}{\pdel x_3},\frac{\pdel}{\pdel x_2}) = 1$, and the
others are zero. ($\ds  x_4=y_1$ and $\ds  x_3= y_2$ 
for the Darboux coordinate 
$\ds  
x_1,x_2,y_1,y_2$
we explained in the section 2). 

We consider the standard basis of homogeneous polynomials of 
$\ds  x_1, x_2, x_3 , x_4$ given by 
$\ds  
\frac{ x_1^{i}}{i!} 
\frac{ x_2^{j}}{j!} 
\frac{ x_3^{k}}{k!} 
\frac{ x_4^{\ell}}{\ell !}$.   
We use the notation 
$\ds  z_{i, j, k, \ell}$ to express the 
dual basis of those.   
While      
we deal with polynomials with $0\le i,j,k,\ell\le 9$, the 4
digit number has the unique meaning and we may denote   
$\ds  z_{i, j, k, \ell}$ by     
$\ds  z_{ijk \ell}$.
Furthermore, in order to simplify
expression we use the notation       
$\ds  \Z{ijk}{i+j+k+\ell}$ for $z_{ijk\ell}$.   

Now, the Poisson bracket is given by 
$$ \Pkt{f}{g} = 
\frac{\pdel f}{\pdel x_1} \frac{\pdel g}{\pdel x_4} 
+\frac{\pdel f}{\pdel x_2} \frac{\pdel g}{\pdel x_3} 
-\frac{\pdel f}{\pdel x_3} \frac{\pdel g}{\pdel x_2} 
-\frac{\pdel f}{\pdel x_4} \frac{\pdel g}{\pdel x_1}.$$ 
\newcommand{\fatte}{\mathbf{e}}
If we denote 
$\ds  
\frac{ x_1^{a_1}}{a_1!} 
\frac{ x_2^{a_2}}{a_2!} 
\frac{ x_3^{a_3}}{a_3!} 
\frac{ x_4^{a_4}}{a_4!}$ by $\ds  \fatte_A$, 
where $\ds  A =(a_1,a_2,a_3,a_4) \in \mZ_{\ge 0}^4$, then 
\begin{align*}
\Pkt{ \fatte_A}{\fatte_B} =& 
+ (a_1b_4-a_4 b_1) 
\frac{ (a_1+b_1-1)!}{a_1 ! b_1 !} 
\frac{ (a_2+b_2  )!}{a_2 ! b_2 !} 
\frac{ (a_3+b_3  )!}{a_3 ! b_3 !} 
\frac{ (a_4+b_4-1)!}{a_4 ! b_4 !} 
\fatte_{(a_1+b_1-1, a_2+b_2, a_3 + b_3, a_4+b_4-1)}\\& 
+ (a_2b_3-a_3 b_2) 
\frac{ (a_1+b_1  )!}{a_1 ! b_1 !} 
\frac{ (a_2+b_2-1)!}{a_2 ! b_2 !} 
\frac{ (a_3+b_3-1)!}{a_3 ! b_3 !} 
\frac{ (a_4+b_4  )!}{a_4 ! b_4 !} 
\fatte_{(a_1+b_1  , a_2+b_2-1, a_3 + b_3-1, a_4+b_4  )}. 
\end{align*}

From the definition of the coboundary operator $\myd$,   
$\ds   \myd z_C ( \fatte_A, \fatte_B ) = - \langle z_C, \Pkt{
\fatte_A}{\fatte_B} \rangle$ holds good, and we see that 
\begin{align*}
\myd z_C = & 
- \sum_{(A,B)\in\text{out}(C)} (a_1 b_4 -a_4 b_1) \frac{C!}{A! B!} z_A \otimes z_B 
- \sum_{(A,B)\in\text{inn}(C)} (a_2 b_3 -a_3 b_2) \frac{C!}{A! B!} z_A \otimes z_B 
\end{align*} 
where $\ds  
C! = c_1 ! c_2 ! c_3 ! c_4 !$,  
$\text{out}(C)$ consists of $(A,B)$ with 
$a_1+b_1 = 1+c_1$, 
$a_2+b_2 =   c_2$, 
$a_3+b_3 =   c_3$, 
$a_4+b_4 = 1+c_4$, 
$|A| > 1$ and 
$|B| > 1$, and also   
$\text{inn}(C)$  consists of $(A,B)$ with 
$a_1+b_1 =   c_1$, 
$a_2+b_2 = 1+c_2$, 
$a_3+b_3 = 1+c_3$, 
$a_4+b_4 =   c_4$, 
$|A| > 1$ and 
$|B| > 1$.
Here the notation 
$\ds  |A|$ means $a_1 + a_2 + a_3 + a_4 $. 
Since we are working in $\frakham_{2n}^{1}$, 
$|A| > 1$ and $|B| > 1${\color{black}.}

\section{Irreducible decomposition of cochain complex}
Our purpose here is to find the multiplicity of the trivial representations for a given
cochain complex. 
We try getting   
a complete irreducible decomposition by finding the maximal weight
vectors.  
A way of finding maximal weight vectors is to
find the invariant vectors by the maximal unipotent subgroup of
$Sp(4,\mR)$.  

\subsection{weight $=2$}  

Since $\ds   \CGFF{1}{4}{1}{2} = \frakS{4}$,  and
 $\ds  \frakS{4}$ is a
non-trivial irreducible representation of $Sp(4,\mR)$, 
$\ds  
\CGF{1}{2} = \{ \mathbf{0} \}$.  
Concerning 
$\ds   \CGFF{2}{4}{1}{2} =  \Lambda^2 \frakS{3}$, 
by finding maximal weight vectors,  
we get an irreducible decomposition 
$$
 \Lambda^2 \frakS{3} =  
\V{0} \oplus \V{1, 1} \oplus \V{2, 2} \oplus \V{3, 3} \oplus
\V{4} \oplus \V{5, 1} \ , $$ 
where $\ds  \V{p,q}$ is the irreducible representation of 
 $Sp(4,\mR)$
corresponding to a Young diagram of length not greater than 2.  
We often denote 
$\ds  \V{p,0}$ by $ \V{p}$, and 
$\ds   \V{p}$ is identical with $\ds  \frakS{p}$,  and 
$\ds  \V{0,0} = \V{0} = 
\frakS{0}$ is the \text{the trivial representation}.    
Hence, 
$\ds  
\CGF{2}{2} \cong  \mR $.  

The following is the reason why we try getting a complete irreducible
decomposition.
By Weyl's dimension formula, we can calculate the dimension
of $\ds  \V{p,q}$ for each $p \ge q \ge 0$ and also
$\ds  \dim \Lambda^r \V{p,q}$,  $\ds  \dim \left(
\Lambda^r \V{p,q} \otimes  \Lambda^{r'} \V{p',q'}\right)$ and so on.  
When we get an irreducible decomposition of $\ds  \Lambda^r
\V{p,q}$, we can calculate the dimension of the left hand side of the
decomposition by using the dimension formula, and also the dimensions of
all terms which appears in the right hand side.  This gives us an
evidence which supports that our decomposition given by a computer
program would be correct.

\begin{kmProp}\label{weight2cochain}
When weight =2, then we have the following table.  
\begin{center}
\begin{tabular}[t]{|c|*{3}{c}|} \hline
degree & 0 & 1 & 2  \\\hline
dim    & 0 & 0 & 1  \\\hline 
\end{tabular}
\end{center} 
Thus, we see the Euler characteristic number of weight 2 is $(-1)^0 0 +
(-1)^1 0 + (-1)^2 1 = 1$.  
\end{kmProp}

Note that we can also see that
in the case of weight $2$,
Betti number $h^2$ of the cohomology is equal to $1$ and $h^0=h^1=0$
by the observation of Proposition~\ref{weight2cochain}.

\subsection{weight = 4 and relative}
The Littlewood-Richardson rule {\color{black}is used} to decompose the tensor product of
two irreducible representations into irreducible components in the
natural representation of $Sp(2n,\mR)$ in general.  
In many places below, 
by the help of Littlewood-Richardson rule, we decompose {\color{black} the} tensor
product $\ds  \V{p,q} \otimes \V{p',q'}$ of two irreducible
representations ({\color{black} of} length at most 2) into irreducible components. We
write down the results only, but by the symbol 
$\myLR$,
we suggest using
of the Littlewood-Richardson rule.  


Here we review the Littlewood-Richardson rule and the specialization
algorithm briefly.  
The product of Schur functions is given by 
$$ s_{\mu}s_{\nu} =\sum_{\lambda} LR_{\mu\nu}^{\lambda} s_{\lambda}\ ,$$ 
where $\lambda, \mu, \nu$ are general partitions and $\displaystyle
LR_{\mu\nu}^{\lambda} $ is called the Littlewood-Richardson coefficient.  
$\displaystyle LR_{\mu\nu}^{\lambda}$ is obtained by
$$ LR_{\mu\nu}^{\lambda} =  \#\{ T\in SSTab(\lambda/\mu: \nu) \mid w(T) \text{ is a lattice permutation}\}$$
(this is called the Littlewood-Richardson rule).  
$SSTab(\lambda/\mu,\nu)$ is the set of semistandard tableaux of shape
$\lambda/\mu$ and weight $\nu$.  
When $T_i$ denotes the $i$-th row of $T$ and $\operatorname{reverse}(T_i)$ the
word obtained by reading $T_i$ from right to left, the $w(T)$ is the
concatenated word
$(\operatorname{reverse}(T_1),\dots,\operatorname{reverse}(T_m))$ whose
length equals $|\nu|=|\lambda|-|\mu|$, where $m=\ell(T)$, the length of $T$.  


In order to obtain the character corresponding to each irreducible
representation of $Sp(2n,\mR)$, the virtual character $\displaystyle
S_{\la\lambda\ra}$ and the universal character $\displaystyle
s_{\la\lambda\ra}$ are defined in the completely same form
as follows.    
%
\begin{align*}
S_{\la\lambda\ra} & = \frac{1}{2}\det( H_{\lambda_i -i+j} +
H_{\lambda_i-i-j})_{1\le i,j\le \ell(\lambda)} \\
s_{\la\lambda\ra} &= \frac{1}{2}\det( h_{\lambda_i -i+j} +
h_{\lambda_i-i-j})_{1\le i,j\le \ell(\lambda)}
\end{align*} 
where $H_k$ is the character of the $k$-th symmetric tensor product of the
natural representation $Sp(2n,\mR)$, $H_k =0$ if $k<0$, and   
$h_k$ is the $k$-th complete symmetric function.  

As the Schur functions $\displaystyle \{ s_{\lambda} \}$ 
is a basis of the space of symmetric functions, 
$\displaystyle \text{Sym}_{\infty}$,  
$\displaystyle \{ s_{\la\lambda\ra}\}$ is also a basis of 
$\displaystyle \text{Sym}_{\infty}$.  %
It is known that 
the product of two universal characters is given by $$
s_{\langle\mu\rangle}s_{\langle\nu\rangle} =
\sum_{\lambda} LR_{\langle\mu\rangle\langle\nu\rangle}^{\langle\lambda\rangle}
s_{\langle\lambda\rangle} \ , 
$$
where 
$$
LR_{\langle\mu\rangle\langle\nu\rangle}^{\langle\lambda\rangle}
= \sum_{\alpha,\beta,\gamma} LR_{\alpha\beta}^{\mu}
LR_{\alpha\gamma}^{\nu} LR_{\beta\gamma}^{\lambda}\ .  $$

There is a homomorphism  
$$\pi: \text{Sym}_{\infty} \rightarrow Rep(Sp(2n,R)) : h_k \mapsto H_k \quad (k=1,2,\ldots) \ ,  
$$ 
which is 
surjective and whose kernel is generated by $e_k-e_{2n-k}$ ($0\le k\le n$) and $e_k$ ($k>2n$), 
where $e_k$ is the $k$-th elementary symmetric function.   
$\pi$ is called the specialization homomorphism and   
satisfies 
$$
\pi ( s_{\langle\lambda\rangle}) = S_{\langle\lambda\rangle} .$$     
In the $Sp(2n,\mR)$-representation theory, each equivalence
classes of the irreducible representations are parametrized by the set 
of partitions $\lambda$ whose length $\ell(\lambda) \le n$.  
We use the notation
$\V{\lambda}$ for that representation.    
Although the character 
$\pi(s_{\langle\lambda\rangle}) = S_{\langle\lambda\rangle}$   
is defined for 
each general partition $\lambda$, if $\ell(\lambda)\le n$ then 
$\pi(s_{\langle\lambda\rangle}) = S_{\langle\lambda\rangle}$   
is the irreducible character of 
$\V{\lambda}$.
\kmcomment{
The symplectic structure yields the isomorphism $\displaystyle
\Lambda^{k} (\mathbb{R}^{2n}) \cong \Lambda^{2n-k}(\mathbb{R}^{2n})$,
and  
}
If $\ell(\lambda) > n$,
it is known $\displaystyle S_{\langle  \lambda\rangle  } = 0 $ or
$\displaystyle \pm S_{\langle  \lambda ' \rangle  }$
with $\ell(\lambda')\le n$ by some rule,
which is called the specialization algorithm.  
%
%
After applying the specialization algorithm completely, we have
$$ S_{\langle \mu \rangle}\otimes S_{\langle \nu \rangle}
=\sum_{\ell(\lambda)\le n} 
LR_{\langle\mu\rangle\langle\nu\rangle}^{\langle\lambda\rangle}S_{\langle \lambda \rangle}
+\sum_{\ell(\lambda)>n} 
LR_{\langle\mu\rangle\langle\nu\rangle}^{\langle\lambda\rangle}S_{\langle \lambda \rangle}
= \sum_{\ell(\lambda)\le n} 
mLR_{\langle\mu\rangle\langle\nu\rangle}^{\langle\lambda\rangle}S_{\langle \lambda \rangle} 
, 
$$
where $mLR_{\langle\mu\rangle\langle\nu\rangle}^{\langle\lambda\rangle}$
is the slight modification of ${LR}_{\langle\mu\rangle\langle\nu\rangle}^{\langle\lambda\rangle}$
caused by applying the specialization algorithm. 
In our context, 
we have 
$$ \V{\mu}\otimes \V{\nu}=\sum_{\ell(\lambda)\le n} 
mLR_{\langle\mu\rangle\langle\nu\rangle}^{\langle\lambda\rangle}\V{\lambda}.$$ 
About the Littlewood-Richardson rule for $Sp(2n,\mR)$, we refer to 
Soichi Okada's book:
{\em Representation theory of classical groups and Combinatorics,
Baifukan, 2006
(in Japanase)} 
Volume2, p.258 and also Volume 2, p.253 for
specialization algorithm.  

We stress that all the 
$ \displaystyle \mathop{\cong}^{\text{LR}} $
in the paper are done by {ourselves} by
following the Littlewood-Richardson rule and the specialization
algorithm faithfully.    

\medskip


The decomposition for weight $=4$ is as follows.  
\begin{align*} 
\CGFF{1}{4}{1}{4} =& \frakS{6}\\ 
\CGFF{2}{4}{1}{4} 
\cong & (\frakS{3}\otimes   
\frakS{5}) \oplus \Lambda^2 \frakS{4}  
\\
\myLR & 
 (\V{2}\oplus\V{4}\oplus\V{6}\oplus\V{8}\oplus\V{3,1}\oplus\V{4,2}
 \oplus\V{5,1}\oplus\V{5,3}\oplus\V{6,2}\oplus\V{7,1}
) \oplus  \Lambda^2 \frakS{4}\\
=&
(\V{2}\oplus\V{4}\oplus\V{6}\oplus
\V{8}\oplus\V{3,1}\oplus\V{4,2}\oplus\V{5,1}\oplus\V{5,3}
\oplus\V{6,2}\oplus\V{7,1}) \\ 
  &\qquad \oplus  
(\V{2}+\V{3,1}+\V{4,2}+\V{5,3}+\V{6}+\V{7,1}
)\\
= & 
2\V{2}+\V{4}+2\V{6}+\V{8}+2\V{3, 1}+2\V{4, 2}+\V{5, 1}+2\V{5, 3}
+\V{6, 2} +2 \V{7, 1} \\
\CGFF{3}{4}{1}{4} \cong &  
 \Lambda^2 \frakS{3}
 \otimes \frakS{4} 
\cong  
( \V{0}+\V{4} + \V{1,1}+\V{2,2}+\V{3,3}+\V{5,1})
 \otimes \V{4} \\
\myLR & + \V{4}\\ 
&+
(
\V{0}+\V{2}+\V{4}+\V{6}+\V{8}+
\V{1,1}+\V{2,2}+\V{3,1}+\V{3,3} 
\\
&\qquad 
+\V{4,2}+\V{4,4}+\V{5,1}+\V{5,3}+\V{6,2}+\V{7,1}
) \\ &+ 
(
\V{4}+\V{3,1}+\V{5,1}
)\\  & +(
\V{4}+\V{2,2}+\V{3,1}+\V{4,2}+\V{5,1}+\V{6,2}
)\\  & +(
\V{4}+\V{3,1}+\V{3,3}+\V{4,2}+\V{5,1}+\V{5,3}+\V{6,2} 
+\V{7,3}
) \\ 
 &+(
\V{2}+\V{4}+\V{6}+\V{8}+\V{1,1}+\V{2,2}+2\V{3,1}+\V{3,3} 
+2\V{4,2}+\V{4,4}\\&\qquad +2\V{5,1}+2\V{5,3}+\V{5,5}+2\V{6,2} 
+\V{6,4}+2\V{7,1}+\V{7,3}+\V{8,2}+\V{9,1}
)
\\
=& 
\V{0}+2\V{2}+6\V{4}+2\V{6}+2\V{8}+2\V{1,1}+3\V{2,2}+6\V{3,1} 
+3\V{3,3}+5\V{4,2}+2\V{4,4}\\&\quad +6\V{5,1}+4\V{5,3}+\V{5,5} 
+5\V{6,2}+\V{6,4}+3\V{7,1}+2\V{7,3}+\V{8,2}+\V{9,1} 
\\
\CGFF{4}{4}{1}{4} = & 
 \Lambda^4 \frakS{3} \\
= & 
3 \V{0} +4 \V{4} +2 \V{6} +\V{8}
+2 \V{1,1} +4 \V{2,2}+3 \V{3,1}+4 \V{3,3} +3 \V{4, 2} \\&
+3 \V{4,4} +5 \V{5,1} +2 \V{5, 3}+\V{5,5} +4 \V{6, 2}
+\V{6, 4}+\V{6, 6}+\V{7,1}+2 \V{7,3} +\V{8,2}.  
\end{align*} 

Thus we obtain 
\begin{align*}
& \CGF{1}{4} = \{ \mathbf{0}\}   
& \CGF{2}{4}& = \{ \mathbf{0}\}\\ 
& \CGF{3}{4} \cong \mR    
& \CGF{4}{4}&\cong \mR^3.\\  
\end{align*}  
\begin{kmProp}
When the weight=4, we have 
\begin{center}
\begin{tabular}{|c|*{5}{c}|} \hline
degree& 0 & 1 & 2 & 3 & 4  \\\hline
dim   & 0 & 0 & 0 & 1 & 3  \\\hline
\end{tabular}
\end{center} 
Thus, the Euler characteristic number of weight 4 is 
$ ( -1)^0 0 + (-1)^3 1 + (-1)^4 3 = 2$.  
\end{kmProp}

\subsection{ weight = 6 and relative} 
The decomposition into irreducible representations for degrees 1,2 and 3
are as follows.  
\begin{align*} 
\CGFF{1}{4}{1}{6} = &  \frakS{8}\\  
\CGFF{2}{4}{1}{6} \cong &
( \frakS{3}\otimes \frakS{7}) \oplus 
( \frakS{4}\otimes \frakS{6}) \oplus 
\Lambda^2 \frakS{5} \\
\myLR & 
(
\V{4}+\V{6}+\V{8}+\V{10}+\V{5, 1}+\V{6, 2}+\V{7, 1}+\V{7, 3}+\V{8, 2}+\V{9, 1}
) \\& 
+
(
\V{2}+\V{4}+\V{6}+\V{8}+\V{10}+\V{3, 1}+\V{4, 2}+\V{5, 1}+\V{5, 3}+\V{6, 2} \\
& \qquad +\V{6, 4}+\V{7, 1}+\V{7, 3}+\V{8, 2}+\V{9, 1}
) \\
& +
(
\V{0}+\V{1,1}+\V{2,2}+\V{3,3}+\V{4}+\V{4,4}+\V{5,1}+\V{5,5}
+\V{6,2}
\\&\qquad 
+\V{7,3}+\V{8}+\V{9,1}
) \\
= & 
\V{0}+\V{2}+3\V{4}+2\V{6}+3\V{8}+2\V{10}+\V{1,1}+\V{2,2}
+\V{3,1}+\V{3,3}+\V{4,2}
\\& 
+\V{4,4}
+3\V{5,1}+\V{5,3}+\V{5,5}
+3\V{6,2}+\V{6,4}+2\V{7,1}+3\V{7,3}+2\V{8,2}+3\V{9,1}\\ 
\CGFF{3}{4}{1}{6} \cong  & 
(
 \Lambda^2 \frakS{3} \otimes \frakS{6}) \oplus (
 \frakS{3} \otimes \frakS{4} \otimes \frakS{5}) \oplus 
 \Lambda^3 \frakS{4} \\
\myLR & 
 (
2 \V{2}+2 \V{4}+6 \V{6}+2 \V{8}+2 \V{10}+\V{1, 1}+\V{2, 2}+3
\V{3, 1}+2 \V{3, 3} \\ 
& \quad 
+5 \V{4, 2}
+\V{4, 4}+6 \V{5, 1}+4 \V{5, 3}+\V{5, 5}+5 \V{6, 2}+3 \V{6, 4}+6
\V{7, 1}\\
&\quad 
+4 \V{7, 3}+\V{7, 5}
+5 \V{8, 2}+\V{8, 4}+3 \V{9, 1}+2 \V{9, 3}+\V{10, 2}+\V{11, 1}
)\\&
+ ( 
\V{0}+6 \V{2}+10 \V{4}+10 \V{6}+6 \V{8}+3 \V{10}+\V{12}+3 \V{1, 1}+5 \V{2, 2}\\
& \quad +12 \V{3, 1}+6 \V{3, 3}+15 \V{4, 2}+5 \V{4, 4}+16 \V{5, 1}+13 \V{5, 3}+3 \V{5, 5} \\ & \quad 
+15 \V{6, 2}+8 \V{6, 4}+\V{6, 6}+12 \V{7, 1}+10 \V{7, 3}+3 \V{7, 5}+9 \V{8, 2}+4 \V{8, 4}\\
&\quad +6 \V{9, 1}+4 \V{9, 3}+3 \V{10, 2}+2 \V{11, 1}
) \\&
+ (
2 \V{2}+3 \V{3, 1}+\V{3, 3}+\V{4}+4 \V{4, 2}+3 \V{5, 1}+3 \V{5, 3}+3 \V{6, 0}+2 \V{6, 2}\\
&\quad +2 \V{6, 4}+2 \V{7, 1}+2 \V{7, 3}+\V{7, 5}+2 \V{8, 2}+\V{9, 1}+\V{9, 3}+\V{10}
) 
\\
= & 
\V{0}+8 \V{2}+12 \V{4}+16 \V{6}+8 \V{8}+5 \V{10}+\V{12}+4 \V{1, 1}+2 \V{2}+6 \V{2, 2}\\
& +18 \V{3, 1}+9 \V{3, 3}+\V{4}+24 \V{4, 2}+6 \V{4, 4}+25 \V{5, 1}+20 \V{5, 3}+4 \V{5, 5}\\
& +3 \V{6}+22 \V{6, 2}+13 \V{6, 4}+\V{6, 6}+20 \V{7, 1}+16 \V{7, 3}+5 \V{7, 5}+16 \V{8, 2}\\
& +5 \V{8, 4}+10 \V{9, 1}+7 \V{9, 3}+\V{10}+4 \V{10, 2}+3 \V{11, 1}\\ 
= & 
\V{0}+10\V{2}+13 \V{4}+19 \V{6}+8 \V{8}+6 \V{10}+\V{12}+4 \V{1, 1} 
+6 \V{2, 2}\\
& +18 \V{3, 1}+9 \V{3, 3}+ 24 \V{4, 2}+6 \V{4, 4}+25 \V{5, 1}+20 \V{5, 3}+4 \V{5, 5}\\
& 
+22 \V{6, 2}+13 \V{6, 4}+\V{6, 6}+20 \V{7, 1}+16 \V{7, 3}+5 \V{7, 5}+16 \V{8, 2}\\
& +5 \V{8, 4}+10 \V{9, 1}+7 \V{9, 3} 
+4 \V{10, 2}+3 \V{11, 1}. 
\end{align*} 
For degree 4, since 
$$ \CGFF{4}{4}{1}{6} 
\cong 
( \Lambda^3 \frakS{3} \otimes \frakS{5}) \oplus 
( \Lambda^2 \frakS{3} \otimes \Lambda^2 \frakS{4}),
$$ 
we shall decompose  two kinds of tensor products:
The first  one is 
\begin{align*}
\Lambda^3 \frakS{3} \otimes \frakS{5} =& 
(\V{2,1}+3 \V{3,0}+\V{3,2}+2 \V{4,1}+\V{4,3}+2 \V{5,2}+\V{6,1} +\V{6,3}+\V{7,0}) \otimes \V{5} 
\\ \myLR & 
9 \V{2}+13 \V{4}+13 \V{6}+10 \V{8}+2 \V{10}+\V{12}+3 \V{1,1}+6
\V{2,2}
+18 \V{3,1}
\\&  +8 \V{3,3}+23 \V{4,2}+7 \V{4,4}+22 \V{5,1}+22 \V{5,3}+4 \V{5,5}+24 \V{6,2}
\\&  
+13 \V{6,4}+2 \V{6,6}+19 \V{7,1}+15 \V{7,3}+6 \V{7,5}+13 \V{8,2}+7 \V{8,4}
\\&  
+\V{8,6}+8 \V{9,1}+7 \V{9,3}+\V{9,5}+5 \V{10,2}+\V{10,4}+2 \V{11,1}+\V{11,3}. 
\end{align*} 
Concerning 
$\ds  \Lambda^2 \frakS{3} \otimes \Lambda^2 \frakS{4}$,   
by using the Littlewood-Richardson rule as usual, 
we have a decomposition into terms
which include a term of partition length 4. 
By an easy observation,  
$\ds  \V{i,j,1,1} = -\V{i,j}$, and we see 
\begin{align*}
\Lambda^2 \frakS{3} \otimes \Lambda^2 \frakS{4} \myLR &   
18 \V{2}+22 \V{4}+26 \V{6}+12 \V{8}+5 \V{10}+\V{12}+6 \V{1,1}+12
\V{2,2}+34 \V{3, 1}\\&  
+15 \V{3,3}+45 \V{4,2}+13 \V{4,4}+41 \V{5,1}+39 \V{5,3}+7
\V{5,5}+40 \V{6,2}\\&  
+25 \V{6 ,4}+3 \V{6,6}+32 \V{7,1}+26 \V{7,3}+11 \V{7,5}+24
\V{8,2}+12 \V{8,4}\\&  
+3 \V{8,6}+12 \V{9 ,1}+12 \V{9,3}+2 \V{9,5}+7 \V{10,2}+3
\V{10,4}+4 \V{11,1}\\&  
+\V{11,3}+\V{12,2}+8 \V{3,2,2 ,1}+\V{3,3,2,2}+4 \V{3,3,3,1}+3
\V{4,2,2,2}\\&  
+12 \V{4,3,2,1}+2 \V{4,3,3,2}+\V{4,4,2,2} +3 \V{4,4,3,1}+8
\V{5,2,2,1}+4 \V{5,3,2,2}\\&  
+4 \V{5,3,3,1}+\V{5,3,3,3}+8 \V{5,4,2,1}+\V{5,4,3,2}+2
\V{5,5,3,1}+\V{6,2,2,2}\\&  
+9 \V{6,3,2,1}+\V{6,3,3,2}+\V{6,4,2,2}+2 \V{6,4,3 ,1}+2
\V{6,5,2,1}+5 \V{7,2,2,1}\\&  
+3 \V{7,3,3,1}+2 \V{7,4,2,1}+2 \V{8,3,2,1}. \\ %
\noalign{
Next, we use 
$\ds  \V{i,j,k,1} = 0$  and 
$\ds  \V{i,j,k,2} = 0$ where $i \ge j \ge k >1 $, and 
have 
} 
\Lambda^2 \frakS{3} \otimes \Lambda^2 \frakS{4} \cong &   
18 \V{2}+22 \V{4}+26 \V{6}+12 \V{8}+5 \V{10}+\V{12}+6 \V{1,1}+12
\V{2,2}\\&  
+34 \V{3, 1}+15 \V{3,3}+45 \V{4,2}+13 \V{4,4}+41 \V{5,1}+39
\V{5,3}+7 \V{5,5}\\&  
+40 \V{6,2}+25 \V{6 ,4}+3 \V{6,6}+32 \V{7,1}+26 \V{7,3}+11
\V{7,5}\\&  
+24 \V{8,2}+12 \V{8,4}+3 \V{8,6}+12 \V{9 ,1}+12 \V{9,3}+2
\V{9,5}\\&  
+7 \V{10,2}+3 \V{10,4}+4
\V{11,1}+\V{11,3}+\V{12,2}+\V{5,3,3,3 }. \\ %
\noalign{
Since $\ds  \V{5,3,3,3} = 0$, we finally have}  
\Lambda^2 \frakS{3} \otimes \Lambda^2 \frakS{4} \cong &   
18 \V{2}+22 \V{4}+26 \V{6}+12 \V{8}+5 \V{10}+\V{12}+6 \V{1,1}+12
\V{2,2}\\&  
+34 \V{3, 1}+15 \V{3,3}+45 \V{4,2}+13 \V{4,4}+41 \V{5,1}+39
\V{5,3}+7 \V{5,5}\\&  
+40 \V{6,2}+25 \V{6 ,4}+3 \V{6,6}+32 \V{7,1}+26 \V{7,3}+11
\V{7,5}\\&  
+24 \V{8,2}+12 \V{8,4}+3 \V{8,6}+12 \V{9 ,1}+12 \V{9,3}+2
\V{9,5}\\&  
+7 \V{10,2}+3 \V{10,4}+4 \V{11,1}+\V{11,3}+\V{12,2}.
\end{align*}
The decomposition is obtained by a computer program.
Validity of our computation is supported  
by the following fact: 
The sum of dimensions of each components above is
113050. On the other hand, 
$\ds  \dim ( \Lambda^2\frakS{3} \otimes\Lambda^2\frakS{4}) =
190\times 595$. These numbers are equal.

Combining above decompositions, we have 
\begin{align*}
\CGFF{4}{4}{1}{6} \cong & 
( \Lambda^3 \frakS{3} \otimes \frakS{5}) \oplus 
( \Lambda^2 \frakS{3} \otimes \Lambda^2 \frakS{4})
\\
=& 
27 \V{2}+35 \V{4}+39 \V{6}+22 \V{8}+7 \V{10}+2 \V{12}+9
\V{1,1}+18 \V{2,2}+52 \V{ 3,1}
\\&   
+23 \V{3,3}+68 \V{4,2}+20 \V{4,4}+63 \V{5,1}+61 \V{5,3}+11 \V{5,5}+64 \V{6,2}
\\&  
+38 \V{6,4}+5 \V{6,6}+51 \V{7,1}+41 \V{7,3}+17 \V{7,5}+37 \V{8,2}+19 \V{8,4}
\\&  
+4 \V{8,6}+20 \V{9,1}+19 \V{9,3}+3 \V{9,5}+12 \V{10,2}+4 \V{10,4} 
+6 \V{11,1}
\\&  
+2 \V{11,3}+\V{12,2}. \\%
\end{align*}

For degree 5, 
applying the Littlewood-Richardson rule to the following decomposition 
\begin{align*} 
& \CGFF{5}{4}{1}{6} \cong  \Lambda^4 \frakS{3} \otimes \frakS{4} 
\cong  
 ( 3 \V{0}+2 \V{1,1}+4 \V{2, 2}+3 \V{3, 1}+4 \V{3, 3}+4 \V{4}
+3 \V{4, 2} 
+3 \V{4, 4}
\\&\quad 
+5 \V{5, 1} +2 \V{5, 3}+\V{5, 5}+2 \V{6}+4 \V{6, 2}+\V{6, 4}
+\V{6, 6}+\V{7, 1}
+2 \V{7, 3}
+\V{8}+\V{8, 2}
) \otimes \V{4} \end{align*} 
and using the Littlewood-Richardson rule many times, we have 
\begin{align*} 
 \V{0} \otimes \V{4} & \cong   \V{4}, \qquad 
 \V{1,1} \otimes \V{4} \myLR   \V{4}+\V{3,1}+\V{5,1}, \\
 \V{2,2} \otimes \V{4}  & \myLR 
  \V{4}+\V{2,2}+\V{3,1}+\V{4,2}+\V{5,1}+\V{6,2}
 \\ 
 \V{3,1} \otimes \V{4} & \myLR 
\V{2}+\V{4}+\V{6}+\V{1,1}+\V{2,2}+2\V{3,1}+\V{3,3}
+2\V{4,2}+2\V{5,1}+\V{5,3}\\&\qquad 
+\V{6,2} +\V{7,1} 
\\ 
 \V{3,3} \otimes \V{4} & \myLR 
  \V{4}+\V{3,1}+\V{3,3}+\V{4,2}+\V{5,1}+\V{5,3}+\V{6,2}
+\V{7,3} 
\\ 
 \V{4} \otimes \V{4} & \myLR 
 \V{0}+\V{2}+\V{4}+\V{6}+\V{8}+\V{1,1}+\V{2,2}+\V{3,1}
+\V{3,3}+\V{4,2}+\V{4,4}\\&\qquad 
+\V{5,1}+\V{5,3}+\V{6,2}+\V{7,1}
\\ 
& \vdots \text{\;\;\;(decompositions of 12 tensor products are omitted)}\\
\kmcomment{
 \V{4,2} \otimes \V{4} & \myLR 
\V{2}+\V{4}+\V{6}+\V{2,2}+2\V{3,1}+\V{3,3}+3\V{4,2}
+\V{4,4}+2\V{5,1}+2\V{5,3}
+2\V{6,2}+\V{6,4}+\V{7,1}+\V{7,3}
+\V{8,2}
\\  
 \V{4,4} \otimes \V{4} & \myLR 
 \V{4}+\V{4,2}+\V{4,4}+\V{5,1}+\V{5,3}+\V{6,2}+\V{6,4}
+\V{7,3}+\V{8,4}
\\ 
 \V{5,1} \otimes \V{4} & \myLR 
 \V{2}+\V{4}+\V{6}+\V{8}+\V{1,1}+\V{2,2}+2\V{3,1}+\V{3,3}
+2\V{4,2}+\V{4,4}+2\V{5,1}
+2\V{5,3}+\V{5,5}+2\V{6,2}
+\V{6,4}+2\V{7,1}+\V{7,3}+\V{8,2}+\V{9,1}
\\ 
 \V{5,3} \otimes \V{4} & \myLR 
 \V{4}+\V{6}+\V{3,1}+\V{3,3}+2\V{4,2}+\V{4,4}+2\V{5,1}
+3\V{5,3}+\V{5,5}+2\V{6,2} 
+2\V{6,4}+\V{7,1}+2\V{7,3}
+\V{7,5}+\V{8,2}+\V{8,4}+\V{9,3}
\\  
\V{5,5} \otimes \V{4} & \myLR 
\V{5,1}+\V{5,3}+\V{5,5}+\V{6,2}+\V{6,4}+\V{7,3}+\V{7,5}
+\V{8,4}+\V{9,5}
\\ 
  \V{6,0} \otimes \V{4} & \myLR 
 \V{2}+\V{4}+\V{6}+\V{8}+\V{10}+\V{3,1}+\V{4,2}+\V{5,1}
+\V{5,3}+\V{6,2}+\V{6,4}+\V{7,1} 
+\V{7,3}+\V{8,2}+\V{9,1}
\\
  \V{6,2} \otimes \V{4} & \myLR 
 \V{4}+\V{6}+\V{8}+\V{2,2}+\V{3,1}+\V{3,3}+2\V{4,2}+\V{4,4}
+2\V{5,1}+2\V{5,3}+\V{5,5} 
+3\V{6,2}+2\V{6,4}+\V{6,6}
+2\V{7,1}+2\V{7,3}+\V{7,5}+2\V{8,2}+\V{8,4}+\V{9,1}+\V{9,3}
+\V{10,2}
\\  
 \V{6,4} \otimes \V{4} & \myLR 
 \V{6}+\V{4,2}+\V{4,4}+\V{5,1}+2\V{5,3}+\V{5,5}+2\V{6,2}
+3\V{6,4}+\V{6,6}+\V{7,1}
+2\V{7,3}+2\V{7,5}+\V{8,2}
+2\V{8,4}+\V{8,6}+\V{9,3}+\V{9,5}+\V{10,4}
\\ 
 \V{6,6} \otimes \V{4} & \myLR 
\V{6,2}+\V{6,4}+\V{6,6}+\V{7,3}+\V{7,5}+\V{8,4}+\V{8,6}
+\V{9,5}+\V{10,6}
\\ 
 \V{7,1} \otimes \V{4} & \myLR 
 \V{4}+\V{6}+\V{8}+\V{10}+\V{3,1}+\V{4,2}+2\V{5,1}+\V{5,3}
+2\V{6,2}+\V{6,4}+2\V{7,1} 
+2\V{7,3}+\V{7,5}+2\V{8,2}
+\V{8,4}+2\V{9,1}+\V{9,3}+\V{10,2}+\V{11,1}
\\ 
  \V{7,3} \otimes \V{4} & \myLR 
 \V{6}+\V{8}+\V{3,3}+\V{4,2}+\V{4,4}+\V{5,1}+2\V{5,3}
+\V{5,5}+2\V{6,2}+2\V{6,4} 
+\V{6,6}+2\V{7,1}+3\V{7,3}
+2\V{7,5}+\V{7,7}+2\V{8,2}+2\V{8,4}+\V{8,6}+\V{9,1} 
+2\V{9,3}+\V{9,5}+\V{10,2}+\V{10,4}+\V{11,3}
\\ 
 \V{8,0} \otimes \V{4} & \myLR 
 \V{4}+\V{6}+\V{8}+\V{10}+\V{12}+\V{5,1}+\V{6,2}+\V{7,1}
+\V{7,3}+\V{8,2}+\V{8,4}+\V{9,1} 
+\V{9,3}+\V{10,2}+\V{11,1}
\\ 
}%
 \V{8,2} \otimes \V{4} & \myLR 
 \V{6}+\V{8}+\V{10}+\V{4,2}+\V{5,1}+\V{5,3}+2\V{6,2}
+\V{6,4}+2\V{7,1}+2\V{7,3} 
\\&\quad 
+\V{7,5} 
+3\V{8,2}+2\V{8,4}
+\V{8,6}+2\V{9,1}+2\V{9,3}+\V{9,5}+2\V{10,2}+\V{10,4}
\\&\quad 
+\V{11,1}+\V{11,3}+\V{12,2}.
\end{align*} 

The complete list of decompositions which we need is in
\ref{Open:two}.

Gathering above decompositions, we see that 
\begin{align*}
\CGFF{5}{4}{1}{6} \cong  & 
4 \V{0}+17 \V{2}+41 \V{4}+29 \V{6}+20 \V{8}+5 \V{10}+\V{12}+12
\V{1, 1}+23 \V{2, 2}\\& 
+45 \V{3, 1}+27 \V{3, 3}+59 \V{4, 2}+24 \V{4, 4}+61 \V{5, 1}+55
\V{5, 3}+15 \V{5, 5}+65 \V{6, 2}\\ & 
+36 \V{6, 4}+8 \V{6, 6}+42 \V{7, 1}+44 \V{7, 3}+16 \V{7, 5}+2 \V{7,7}+31 \V{8, 2}+21 \V{8, 4}\\&
+5 \V{8, 6}+18 \V{9, 1}+15 \V{9, 3}+6 \V{9, 5}+10 \V{10, 2}+4
\V{10, 4}+\V{10, 6}+3 \V{11, 1}\\&
+3 \V{11, 3}+\V{12, 2}. \\[2mm]
\noalign{For degree 6, direct computation of maximal weight vectors shows
that}\\[-3mm]
\CGFF{6}{4}{1}{6} \cong & 
 \Lambda^6 \frakS{3}\\
\cong & 
 4 \V{0}  
+ 6 \V{1,1}  
+ 2 \V{2,0}  
+ 10 \V{2,2}  
+ 10 \V{3,1}  
+ 12 \V{3,3}  
+ 13 \V{4,0}  
+ 14 \V{4,2} \\&  
+  9 \V{4,4}  
+ 19 \V{5,1}  
+ 14 \V{5,3}  
+  7 \V{5,5}  
+  7 \V{6,0}  
+ 18 \V{6,2}  
+  9 \V{6,4}  
+  4 \V{6,6} \\& 
+ 10 \V{7,1}  
+ 13 \V{7,3}  
+  4 \V{7,5}  
+  2 \V{7,7}  
+  4 \V{8,0}  
+  7 \V{8,2}  
+  5 \V{8,4}  
+    \V{8,6} \\&     
+  3 \V{9,1}
+  3 \V{9,3}
+  2 \V{9,5}  
+    \V{10,0}  
+    \V{10,2}.  
\end{align*}  
From the observation above, we have the following:
\begin{kmProp}
When the weight=6, we have

\begin{center}
\begin{tabular}{|c|*{7}{c}|} \hline
degree& 0 & 1 & 2 & 3 & 4 & 5 & 6 \\\hline
dim &   0 & 0 & 1 & 1 & 0 & 4 & 4 \\\hline
\end{tabular}
\end{center} 

The Euler characteristic number for weight 6 is 
$\ds  (-1)^0 0 + (-1)^2 1 +  (-1)^3 1 +  
(-1)^5 4 + 
(-1)^6 4 =0 $.   
\end{kmProp}

\section{Betti numbers}
In order to get each Betti number, we have to know the image and the
kernel of $\myd$ itself.  For that purpose, we have to fix some bases of
cochain complexes and matrix representation of $\myd$ by those concrete
bases.

\subsection{weight=2}

We already see that the Betti numbers of the cohomology group in the
case of weight $2$ are the following which is an immediate consequence of 
Proposition~\ref{weight2cochain}.

\begin{center}
\begin{tabular}[t]{|c|*{3}{c}|}
\hline
degree & 0 & 1 & 2  \\\hline
dim    & 0 & 0 & 1  \\\hline
Betti  & 0 & 0 & 1  \\\hline
\end{tabular}
\end{center}

\subsection{weight=4} 
In order to know properties of 
$\ds  \myd : 
\CGF{3}{4} \rightarrow 
\CGF{4}{4} 
$,
we prepare a basis, say $A$  of 
$\ds  
\CGF{3}{4}$: 
\begin{align*} A = &  -4 \Z{1 0 0}{3} \wedge \Z{2 1 0}{3} \wedge \Z{1 0 1}{4}
+2 \Z{0 0 1}{3} \wedge \Z{1 0 0}{3} \wedge \Z{3 1 0}{4}
+4 \Z{1 0 0}{3} \wedge \Z{2 0 1}{3} \wedge \Z{1 1 0}{4} 
-2 \Z{0 0 1}{3} \wedge \Z{0 2 1}{3} \wedge \Z{2 1 1}{4}
\\&
+\Z{0 0 1}{3} \wedge \Z{0 3 0}{3} \wedge \Z{2 0 2}{4}
-2 \Z{1 0 0}{3} \wedge \Z{3 0 0}{3} \wedge \Z{1 0 0}{4}
-4 \Z{0 0 1}{3} \wedge \Z{2 0 1}{3} \wedge \Z{1 2 0}{4}
+4 \Z{0 0 1}{3} \wedge \Z{2 1 0}{3} \wedge \Z{1 1 1}{4}
\\&
\pm \textbf{164 terms} 
\\& 
-2 \Z{0 0 0}{3} \wedge \Z{2 0 1}{3} \wedge \Z{2 1 0}{4}
+2 \Z{0 0 0}{3} \wedge \Z{2 1 0}{3} \wedge \Z{2 0 1}{4}
+\Z{0 0 0}{3} \wedge \Z{3 0 0}{3} \wedge \Z{2 0 0}{4} 
+\Z{0 0 1}{3} \wedge \Z{0 1 0}{3} \wedge \Z{4 0 0}{4}
\\&
+2 \Z{0 0 1}{3} \wedge \Z{0 1 1}{3} \wedge \Z{3 1 0}{4}
+\Z{0 0 1}{3} \wedge \Z{0 1 2}{3} \wedge \Z{2 2 0}{4} 
\end{align*}

The complete expression of all terms in $A$ is exhibited in 
\ref{Open:wt:four}.

We prepare a basis $\{B_1, B_2, B_3\}$ of $\ds  
\CGF{4}{4}$ and we get $\ds   \myd A = -32 B_1 + 32 B_2 -28 B_3 $
by using computer programs.
The complete form of its basis is also exhibited in 
\ref{Open:wt:four}.  
But we can conclude $\myd A\neq 0$ directly by an observation of 
the calculation of the $\myd$-image of $A$.
Thus,  
$\ds  h^{3}=0$ and $\ds  h^{4}=2$,
 that is, $\ds  h^{0}= h^{1}=h^{2}=h^{3} =0$.
Consequently, $\ds   h^{4}=2$ and the others are zero.
Therefore , the alternating sum of the Betti
numbers, which is another definition of the Euler characteristic number, is 
2.  

\begin{thm}
When the weight=4, we have 
\begin{center}
\begin{tabular}{|c|*{5}{c}|}
\hline
degree& 0 & 1 & 2 & 3 & 4  \\\hline
dim   & 0 & 0 & 0 & 1 & 3  \\\hline
Betti & 0 & 0 & 0 & 0 & 2  \\\hline
\end{tabular}
\end{center}

\end{thm}

\subsection{weight=6}

As a basis of $\ds  \CGF{2}{6}$, we can select  
\begin{align*} A =  &
6 \Z{112}{5}\wedge \Z{121}{5}
+\frac{1}{10} \Z{000}{5}\wedge \Z{500}{5}
+\frac{1}{2} \Z{004}{5}\wedge \Z{140}{5} 
+\frac{1}{2} \Z{001}{5}\wedge \Z{410}{5}  
+\Z{002}{5}\wedge \Z{320}{5}
+\Z{003}{5}\wedge \Z{230}{5}\\&
-2 \Z{011}{5}\wedge \Z{311}{5}
+\frac{1}{10} \Z{005}{5}\wedge \Z{050}{5} 
-\frac{1}{2} \Z{010}{5}\wedge \Z{401}{5}
-3 \Z{012}{5}\wedge \Z{221}{5}
-2 \Z{013}{5}\wedge \Z{131}{5}
-\frac{1}{2} \Z{014}{5}\wedge \Z{041}{5}\\&
+\Z{020}{5}\wedge \Z{302}{5}
+3 \Z{021}{5}\wedge \Z{212}{5} 
+3 \Z{022}{5}\wedge \Z{122}{5} 
+\Z{023}{5}\wedge \Z{032}{5} 
-\Z{030}{5}\wedge \Z{203}{5}
-2 \Z{031}{5}\wedge \Z{113}{5}\\&
+\frac{1}{2} \Z{040}{5}\wedge \Z{104}{5}
-\frac{1}{2} \Z{100}{5}\wedge \Z{400}{5} 
-2 \Z{101}{5}\wedge \Z{310}{5}
-3 \Z{102}{5}\wedge \Z{220}{5}
-2 \Z{103}{5}\wedge \Z{130}{5}\\& 
+2 \Z{110}{5}\wedge \Z{301}{5}
+6 \Z{111}{5}\wedge \Z{211}{5}
-3 \Z{120}{5}\wedge \Z{202}{5}
+\Z{200}{5}\wedge \Z{300}{5}
+3 \Z{201}{5}\wedge \Z{210}{5} .
\end{align*} 

We can also see $\ds  \myd A \neq 0$ by an observation for $A$. 
Thus, $\ds  h^{2} =0$ and $\ds  h^{3} =0$.
We also prepared a basis $B$ of $\ds  \CGF{3}{6}$ and we got $\ds  \myd A = 12 B$,
the complete expression of $B$ is 
in \ref{Open:wt:six}.

\medskip

Concerning the coboundary operator $\ds  \myd : \CGF{5}{6}
\rightarrow  \CGF{6}{6}$, 
we have concrete bases 
$\ds  P_{1}, P_{2}, P_{3}, P_{4}$ of 
 $\ds  \CGF{5}{6}$
and 
$\ds  Q_{1}, Q_{2}, Q_{3}, Q_{4}$ of  
$\ds  \CGF{6}{6}$.
Those have very long expressions as 
$\ds  P_{1}$ is a sum of 3696 terms,   
$\ds  P_{2}$ of 3358 terms,  
$\ds  P_{3}$ of 1406 terms,    
$\ds  P_{4}$ of 2960 terms, and  
$\ds  Q_{1}$ of 120 terms,   
$\ds  Q_{2}$ of 466 terms,  
$\ds  Q_{3}$ of 756 terms,    
$\ds  Q_{4}$ of 866 terms.  
Each $P_i$ is constructed by the terms of the form $\left(\wedge_{a=1}^4 \Z{i_a j_a k_a}{3}\right)\wedge \Z{i_5 j_5 k_5}{4}$,
and each $Q_i$ is by the terms of $\wedge_{a=1}^6 \Z{i_a j_a k_a}{3}$.
The complete forms of each $P_i$ and $Q_i$ are also exhibited in 
\ref{Open:wt:six}. 

Using those bases, 
we have  
a matrix representation of $\myd${\color{black}:} 
\begin{align*}
\myd( P_{1} ) & = -90 Q_{1} -22 Q_{2} + 5 Q_{3} -2 Q_{4} \\ 
\myd( P_{2} ) & = -\frac{25}{2} Q_{1} - \frac{7}{6}  Q_{2} +
\frac{55}{12} Q_{3} +3 Q_{4} \\ 
\myd( P_{3} ) & = -\frac{17}{2} Q_{1} - \frac{4}{3}  Q_{2} +
\frac{7}{6} Q_{3}  \\ 
\myd( P_{4} ) & = 6 Q_{1} + \frac{3}{2}  Q_{2} -
\frac{9}{2} Q_{3} -4 Q_{4}{\color{black}.} 
\end{align*} 
It is non-singular and so 
$\ds  h^{5} =0$ and $\ds  h^{6} =0$. Namely,   
$\ds  h^{j} =0 $ for $j=0,..,6$.

\begin{thm}
When the weight=6, we have

\begin{center}
\begin{tabular}{|c|*{7}{c}|}
\hline
degree& 0 & 1 & 2 & 3 & 4 & 5 & 6 \\\hline
dim &   0 & 0 & 1 & 1 & 0 & 4 & 4 \\\hline
Betti&  0 & 0 & 0 & 0 & 0 & 0 & 0 \\\hline
\end{tabular}
\end{center}
\end{thm}

\nocite{metoki:shinya} 
\nocite{morita:text}
\nocite{morita:text:eng}
\nocite{okada:text}
\nocite{M:Takamura}
\nocite{fulton:harris}
\nocite{goodman:wallach}
\nocite{MR0425981}
\nocite{MR0266195}

\bibliographystyle{plain}
\bibliography{km_refs}

\appendix
\def\thesection{Appendix \Alph{section}}
\section{The complete list needed in subsection 4.3}\label{Open:two}
The list below consists of some decompositions, which we need.     
\begin{align*} 
 \V{0} \otimes \V{4}  \cong &  \V{4}, \qquad 
 \V{1,1} \otimes \V{4} \myLR   \V{4}+\V{3,1}+\V{5,1}, \\
 \V{2,2} \otimes \V{4}   \myLR &
  \V{4}+\V{2,2}+\V{3,1}+\V{4,2}+\V{5,1}+\V{6,2}
 \\ 
 \V{3,1} \otimes \V{4}  \myLR &
\V{2}+\V{4}+\V{6}+\V{1,1}+\V{2,2}+2\V{3,1}+\V{3,3}
+2\V{4,2}+2\V{5,1}+\V{5,3}\\&
+\V{6,2} +\V{7,1} 
\\ 
 \V{3,3} \otimes \V{4}  \myLR &
  \V{4}+\V{3,1}+\V{3,3}+\V{4,2}+\V{5,1}+\V{5,3}+\V{6,2}
+\V{7,3} 
\\ 
 \V{4} \otimes \V{4}  \myLR &
 \V{0}+\V{2}+\V{4}+\V{6}+\V{8}+\V{1,1}+\V{2,2}+\V{3,1}
+\V{3,3}+\V{4,2}+\V{4,4}\\
&
+\V{5,1}+\V{5,3}+\V{6,2}+\V{7,1}
\\ 
 \V{4,2} \otimes \V{4}  \myLR &
\V{2}+\V{4}+\V{6}+\V{2,2}+2\V{3,1}+\V{3,3}+3\V{4,2}
+\V{4,4}+2\V{5,1}+2\V{5,3}\\&
+2\V{6,2}+\V{6,4}+\V{7,1}+\V{7,3}
+\V{8,2}
\\  
 \V{4,4} \otimes \V{4}  \myLR &
 \V{4}+\V{4,2}+\V{4,4}+\V{5,1}+\V{5,3}+\V{6,2}+\V{6,4}
+\V{7,3}+\V{8,4}
\\ 
 \V{5,1} \otimes \V{4}  \myLR &
 \V{2}+\V{4}+\V{6}+\V{8}+\V{1,1}+\V{2,2}+2\V{3,1}+\V{3,3}
+2\V{4,2}+\V{4,4}\\& +2\V{5,1}
+2\V{5,3}
+\V{5,5}+2\V{6,2}
+\V{6,4}+2\V{7,1}+\V{7,3}+\V{8,2}+\V{9,1}
\\ 
 \V{5,3} \otimes \V{4}  \myLR &
 \V{4}+\V{6}+\V{3,1}+\V{3,3}+2\V{4,2}+\V{4,4}+2\V{5,1}
+3\V{5,3}+\V{5,5}+2\V{6,2} 
+2\V{6,4}
\\&
+\V{7,1}+2\V{7,3}
+\V{7,5}+\V{8,2}+\V{8,4}+\V{9,3}
\\  
\V{5,5} \otimes \V{4}  \myLR &
\V{5,1}+\V{5,3}+\V{5,5}+\V{6,2}+\V{6,4}+\V{7,3}+\V{7,5}
+\V{8,4}+\V{9,5}
\\ 
  \V{6,0} \otimes \V{4}  \myLR &
 \V{2}+\V{4}+\V{6}+\V{8}+\V{10}+\V{3,1}+\V{4,2}+\V{5,1}
+\V{5,3}+\V{6,2}+\V{6,4}\\&
+\V{7,1} 
+\V{7,3}+\V{8,2}+\V{9,1}
\\
  \V{6,2} \otimes \V{4}  \myLR &
 \V{4}+\V{6}+\V{8}+\V{2,2}+\V{3,1}+\V{3,3}+2\V{4,2}+\V{4,4}
+2\V{5,1}+2\V{5,3}+\V{5,5} 
+3\V{6,2}\\&
+2\V{6,4}+\V{6,6}
+2\V{7,1}+2\V{7,3}+\V{7,5}+2\V{8,2}+\V{8,4}+\V{9,1}+\V{9,3}
+\V{10,2}
\\  
 \V{6,4} \otimes \V{4}  \myLR &
 \V{6}+\V{4,2}+\V{4,4}+\V{5,1}+2\V{5,3}+\V{5,5}+2\V{6,2}
+3\V{6,4}+\V{6,6}+\V{7,1}
\\&
+2\V{7,3}+2\V{7,5}+\V{8,2}
+2\V{8,4}+\V{8,6}+\V{9,3}+\V{9,5}+\V{10,4}
\\ 
 \V{6,6} \otimes \V{4}  \myLR &
\V{6,2}+\V{6,4}+\V{6,6}+\V{7,3}+\V{7,5}+\V{8,4}+\V{8,6}
+\V{9,5}+\V{10,6}
\\ 
 \V{7,1} \otimes \V{4}  \myLR &
 \V{4}+\V{6}+\V{8}+\V{10}+\V{3,1}+\V{4,2}+2\V{5,1}+\V{5,3}
+2\V{6,2}+\V{6,4}+2\V{7,1} 
\\&
+2\V{7,3}+\V{7,5}+2\V{8,2}
+\V{8,4}+2\V{9,1}+\V{9,3}+\V{10,2}+\V{11,1}
\\ 
  \V{7,3} \otimes \V{4}  \myLR &
 \V{6}+\V{8}+\V{3,3}+\V{4,2}+\V{4,4}+\V{5,1}+2\V{5,3}
+\V{5,5}+2\V{6,2}+2\V{6,4} 
+\V{6,6}\\&
+2\V{7,1}+3\V{7,3}
+2\V{7,5}+\V{7,7}+2\V{8,2}+2\V{8,4}+\V{8,6}+\V{9,1} 
\\&
+2\V{9,3}+\V{9,5}+\V{10,2}+\V{10,4}+\V{11,3}
\\ 
 \V{8,0} \otimes \V{4}  \myLR &
 \V{4}+\V{6}+\V{8}+\V{10}+\V{12}+\V{5,1}+\V{6,2}+\V{7,1}
+\V{7,3}+\V{8,2}+\V{8,4}+\V{9,1} 
\\&
+\V{9,3}+\V{10,2}+\V{11,1}
\\ 
 \V{8,2} \otimes \V{4}  \myLR &
 \V{6}+\V{8}+\V{10}+\V{4,2}+\V{5,1}+\V{5,3}+2\V{6,2}
+\V{6,4}+2\V{7,1}+2\V{7,3} 
\\&\quad 
+\V{7,5} 
+3\V{8,2}+2\V{8,4}
+\V{8,6}+2\V{9,1}+2\V{9,3}+\V{9,5}+2\V{10,2}+\V{10,4}
\\&\quad 
+\V{11,1}+\V{11,3}+\V{12,2}
\end{align*}

\def\thesection{Appendix \Alph{section}}
\section{Concrete bases for weight $=4$}\label{Open:wt:four}
In order to get a matrix representation of 
$\ds\myd : 
\CGF{3}{4} \rightarrow \CGF{4}{4} $: 
we prepare a basis, say $A$  of 
$\ds \CGF{3}{4}$: 
\begin{align*} & A := \\  
&
-4 z^{3}_{1 0 0} \wedge z^{3}_{2 1 0} \wedge z^{4}_{1 0 1}
+2 z^{3}_{0 0 1} \wedge z^{3}_{1 0 0} \wedge z^{4}_{3 1 0}
+4 z^{3}_{1 0 0} \wedge z^{3}_{2 0 1} \wedge z^{4}_{1 1 0} 
-2 z^{3}_{0 0 1} \wedge z^{3}_{0 2 1} \wedge z^{4}_{2 1 1}
+z^{3}_{0 0 1} \wedge z^{3}_{0 3 0} \wedge z^{4}_{2 0 2}\\&
-2 z^{3}_{1 0 0} \wedge z^{3}_{3 0 0} \wedge z^{4}_{1 0 0}
-4 z^{3}_{0 0 1} \wedge z^{3}_{2 0 1} \wedge z^{4}_{1 2 0}
+4 z^{3}_{0 0 1} \wedge z^{3}_{2 1 0} \wedge z^{4}_{1 1 1}
-2 z^{3}_{0 0 1} \wedge z^{3}_{0 2 0} \wedge z^{4}_{3 0 1}
-4 z^{3}_{0 0 1} \wedge z^{3}_{1 1 1} \wedge z^{4}_{1 2 1} \\&
-2 z^{3}_{0 0 1} \wedge z^{3}_{1 1 1} \wedge z^{4}_{2 1 0}
+2 z^{3}_{0 0 1} \wedge z^{3}_{1 2 0} \wedge z^{4}_{1 1 2}
+2 z^{3}_{0 0 1} \wedge z^{3}_{1 2 0} \wedge z^{4}_{2 0 1}
-4 z^{3}_{0 0 1} \wedge z^{3}_{2 0 0} \wedge z^{4}_{2 1 0}
+4 z^{3}_{0 0 2} \wedge z^{3}_{1 2 0} \wedge z^{4}_{1 1 1} \\&
-2 z^{3}_{0 0 2} \wedge z^{3}_{2 0 0} \wedge z^{4}_{1 2 0}
+4 z^{3}_{0 0 1} \wedge z^{3}_{1 0 1} \wedge z^{4}_{2 2 0}
+2 z^{3}_{0 0 1} \wedge z^{3}_{1 0 2} \wedge z^{4}_{1 3 0}
-4 z^{3}_{0 0 1} \wedge z^{3}_{1 1 0} \wedge z^{4}_{2 1 1}
-2 z^{3}_{0 0 1} \wedge z^{3}_{1 1 0} \wedge z^{4}_{3 0 0} \\&
+4 z^{3}_{0 0 2} \wedge z^{3}_{0 1 1} \wedge z^{4}_{2 2 0}
+z^{3}_{0 0 2} \wedge z^{3}_{1 0 0} \wedge z^{4}_{2 2 0}
+2 z^{3}_{0 0 2} \wedge z^{3}_{1 0 1} \wedge z^{4}_{1 3 0}
+z^{3}_{0 0 2} \wedge z^{3}_{1 0 2} \wedge z^{4}_{0 4 0}
-2 z^{3}_{0 0 2} \wedge z^{3}_{1 1 0} \wedge z^{4}_{1 2 1} \\&
-4 z^{3}_{0 0 2} \wedge z^{3}_{1 1 0} \wedge z^{4}_{2 1 0}
-2 z^{3}_{0 0 2} \wedge z^{3}_{1 1 1} \wedge z^{4}_{0 3 1}
-4 z^{3}_{0 0 2} \wedge z^{3}_{1 1 1} \wedge z^{4}_{1 2 0}
+z^{3}_{0 0 2} \wedge z^{3}_{1 2 0} \wedge z^{4}_{0 2 2}
+2 z^{3}_{0 1 0} \wedge z^{3}_{0 1 2} \wedge z^{4}_{2 1 1} \\&
-z^{3}_{0 1 0} \wedge z^{3}_{0 2 1} \wedge z^{4}_{2 0 2}
-2 z^{3}_{0 1 0} \wedge z^{3}_{1 0 2} \wedge z^{4}_{1 2 1}
+2 z^{3}_{0 1 0} \wedge z^{3}_{1 0 2} \wedge z^{4}_{2 1 0}
+4 z^{3}_{0 1 0} \wedge z^{3}_{1 1 0} \wedge z^{4}_{2 0 2}
+z^{3}_{0 0 1} \wedge z^{3}_{2 1 0} \wedge z^{4}_{2 0 0} \\&
+2 z^{3}_{0 0 1} \wedge z^{3}_{3 0 0} \wedge z^{4}_{1 1 0}
+2 z^{3}_{0 0 2} \wedge z^{3}_{0 1 0} \wedge z^{4}_{3 1 0}
+2 z^{3}_{0 0 2} \wedge z^{3}_{0 1 2} \wedge z^{4}_{1 3 0}
-4 z^{3}_{0 0 2} \wedge z^{3}_{0 2 0} \wedge z^{4}_{2 1 1}
-4 z^{3}_{0 0 2} \wedge z^{3}_{0 2 1} \wedge z^{4}_{1 2 1} \\& 
+2 z^{3}_{0 0 2} \wedge z^{3}_{0 3 0} \wedge z^{4}_{1 1 2}
+4 z^{3}_{0 1 1} \wedge z^{3}_{2 0 0} \wedge z^{4}_{1 1 1}
-2 z^{3}_{1 2 0} \wedge z^{3}_{2 0 1} \wedge z^{4}_{0 0 1}
-2 z^{3}_{2 0 0} \wedge z^{3}_{2 0 1} \wedge z^{4}_{0 1 0}
+2 z^{3}_{2 0 0} \wedge z^{3}_{2 1 0} \wedge z^{4}_{0 0 1} \\&
+z^{3}_{2 0 0} \wedge z^{3}_{3 0 0} \wedge z^{4}_{0 0 0}
+z^{3}_{2 0 1} \wedge z^{3}_{2 1 0} \wedge z^{4}_{0 0 0}
-4 z^{3}_{1 1 0} \wedge z^{3}_{1 1 1} \wedge z^{4}_{0 1 2}
+4 z^{3}_{1 1 0} \wedge z^{3}_{1 1 1} \wedge z^{4}_{1 0 1}
+2 z^{3}_{1 1 0} \wedge z^{3}_{1 2 0} \wedge z^{4}_{0 0 3} \\&
-2 z^{3}_{1 1 0} \wedge z^{3}_{2 0 0} \wedge z^{4}_{1 0 1}
-4 z^{3}_{1 1 0} \wedge z^{3}_{2 0 1} \wedge z^{4}_{0 1 1}
+2 z^{3}_{1 1 0} \wedge z^{3}_{2 0 1} \wedge z^{4}_{1 0 0}
+4 z^{3}_{1 1 0} \wedge z^{3}_{2 1 0} \wedge z^{4}_{0 0 2}
+2 z^{3}_{1 1 1} \wedge z^{3}_{2 0 0} \wedge z^{4}_{0 1 1} \\&
+2 z^{3}_{1 1 1} \wedge z^{3}_{2 0 1} \wedge z^{4}_{0 1 0}
+2 z^{3}_{1 1 0} \wedge z^{3}_{3 0 0} \wedge z^{4}_{0 0 1}
+4 z^{3}_{1 1 1} \wedge z^{3}_{1 2 0} \wedge z^{4}_{0 0 2}
+2 z^{3}_{1 1 1} \wedge z^{3}_{2 1 0} \wedge z^{4}_{0 0 1}
-z^{3}_{1 2 0} \wedge z^{3}_{2 0 0} \wedge z^{4}_{0 0 2} \\&
+4 z^{3}_{1 0 2} \wedge z^{3}_{1 1 0} \wedge z^{4}_{1 1 0}
+4 z^{3}_{1 0 2} \wedge z^{3}_{1 1 1} \wedge z^{4}_{0 2 0}
-4 z^{3}_{1 0 2} \wedge z^{3}_{1 2 0} \wedge z^{4}_{0 1 1}
-z^{3}_{1 0 2} \wedge z^{3}_{2 0 0} \wedge z^{4}_{0 2 0}
-2 z^{3}_{1 0 2} \wedge z^{3}_{2 1 0} \wedge z^{4}_{0 1 0} \\&
+4 z^{3}_{1 0 1} \wedge z^{3}_{2 0 1} \wedge z^{4}_{0 2 0}
-4 z^{3}_{1 0 1} \wedge z^{3}_{2 1 0} \wedge z^{4}_{0 1 1}
-2 z^{3}_{1 0 1} \wedge z^{3}_{2 1 0} \wedge z^{4}_{1 0 0}
-2 z^{3}_{1 0 1} \wedge z^{3}_{3 0 0} \wedge z^{4}_{0 1 0}
-2 z^{3}_{1 0 2} \wedge z^{3}_{1 1 0} \wedge z^{4}_{0 2 1} \\&
+3 z^{3}_{1 0 0} \wedge z^{3}_{2 0 0} \wedge z^{4}_{2 0 0}
+4 z^{3}_{1 0 1} \wedge z^{3}_{1 1 0} \wedge z^{4}_{2 0 0}
+4 z^{3}_{1 0 1} \wedge z^{3}_{1 1 1} \wedge z^{4}_{0 2 1}
+4 z^{3}_{1 0 1} \wedge z^{3}_{1 1 1} \wedge z^{4}_{1 1 0}
-2 z^{3}_{1 0 1} \wedge z^{3}_{1 2 0} \wedge z^{4}_{0 1 2} \\&
-4 z^{3}_{1 0 1} \wedge z^{3}_{1 2 0} \wedge z^{4}_{1 0 1}
+2 z^{3}_{1 0 1} \wedge z^{3}_{2 0 0} \wedge z^{4}_{1 1 0}
-2 z^{3}_{0 3 0} \wedge z^{3}_{1 1 1} \wedge z^{4}_{0 0 3}
-2 z^{3}_{1 0 0} \wedge z^{3}_{1 0 1} \wedge z^{4}_{2 1 0}
-z^{3}_{0 3 0} \wedge z^{3}_{2 0 1} \wedge z^{4}_{0 0 2} \\&
-2 z^{3}_{1 0 0} \wedge z^{3}_{1 0 2} \wedge z^{4}_{1 2 0}
+2 z^{3}_{1 0 0} \wedge z^{3}_{1 1 0} \wedge z^{4}_{2 0 1}
+4 z^{3}_{1 0 0} \wedge z^{3}_{1 1 1} \wedge z^{4}_{1 1 1}
-2 z^{3}_{1 0 0} \wedge z^{3}_{1 2 0} \wedge z^{4}_{1 0 2}
+2 z^{3}_{0 2 0} \wedge z^{3}_{0 2 1} \wedge z^{4}_{1 0 3} \\&
+z^{3}_{0 2 0} \wedge z^{3}_{1 0 2} \wedge z^{4}_{0 2 2}
-4 z^{3}_{0 2 0} \wedge z^{3}_{1 0 2} \wedge z^{4}_{1 1 1}
-2 z^{3}_{0 2 0} \wedge z^{3}_{1 1 0} \wedge z^{4}_{1 0 3}
-2 z^{3}_{0 2 0} \wedge z^{3}_{1 1 1} \wedge z^{4}_{0 1 3}
+4 z^{3}_{0 2 0} \wedge z^{3}_{1 1 1} \wedge z^{4}_{1 0 2} \\&
+z^{3}_{0 2 0} \wedge z^{3}_{1 2 0} \wedge z^{4}_{0 0 4}
+2 z^{3}_{0 2 0} \wedge z^{3}_{1 0 1} \wedge z^{4}_{1 1 2}
-4 z^{3}_{0 2 0} \wedge z^{3}_{1 0 1} \wedge z^{4}_{2 0 1}
-2 z^{3}_{0 2 0} \wedge z^{3}_{2 0 0} \wedge z^{4}_{1 0 2}
-2 z^{3}_{0 2 0} \wedge z^{3}_{2 0 1} \wedge z^{4}_{0 1 2} \\&
+2 z^{3}_{0 2 0} \wedge z^{3}_{2 0 1} \wedge z^{4}_{1 0 1}
+2 z^{3}_{0 2 0} \wedge z^{3}_{2 1 0} \wedge z^{4}_{0 0 3}
+z^{3}_{0 2 0} \wedge z^{3}_{3 0 0} \wedge z^{4}_{0 0 2}
+z^{3}_{0 2 1} \wedge z^{3}_{0 3 0} \wedge z^{4}_{0 0 4}
-4 z^{3}_{0 2 1} \wedge z^{3}_{1 0 1} \wedge z^{4}_{1 1 1} \\&
-4 z^{3}_{0 2 1} \wedge z^{3}_{1 0 2} \wedge z^{4}_{0 2 1}
-2 z^{3}_{0 2 1} \wedge z^{3}_{1 1 0} \wedge z^{4}_{1 0 2}
+2 z^{3}_{0 2 1} \wedge z^{3}_{1 1 1} \wedge z^{4}_{0 1 2}
+2 z^{3}_{0 2 1} \wedge z^{3}_{1 2 0} \wedge z^{4}_{0 0 3}
+2 z^{3}_{0 2 1} \wedge z^{3}_{2 0 1} \wedge z^{4}_{0 1 1} \\&
+z^{3}_{0 2 1} \wedge z^{3}_{2 1 0} \wedge z^{4}_{0 0 2}
+2 z^{3}_{0 3 0} \wedge z^{3}_{1 0 1} \wedge z^{4}_{1 0 2}
+2 z^{3}_{0 3 0} \wedge z^{3}_{1 0 2} \wedge z^{4}_{0 1 2}
+2 z^{3}_{0 1 2} \wedge z^{3}_{1 0 1} \wedge z^{4}_{1 2 0}
+2 z^{3}_{0 1 2} \wedge z^{3}_{1 0 2} \wedge z^{4}_{0 3 0} \\&
+4 z^{3}_{0 1 2} \wedge z^{3}_{1 1 0} \wedge z^{4}_{1 1 1}
+2 z^{3}_{0 1 2} \wedge z^{3}_{1 1 1} \wedge z^{4}_{0 2 1}
-4 z^{3}_{0 1 2} \wedge z^{3}_{1 2 0} \wedge z^{4}_{0 1 2}
-z^{3}_{0 1 2} \wedge z^{3}_{2 0 1} \wedge z^{4}_{0 2 0}
-2 z^{3}_{0 1 2} \wedge z^{3}_{2 1 0} \wedge z^{4}_{0 1 1} \\&
+z^{3}_{0 2 0} \wedge z^{3}_{1 0 0} \wedge z^{4}_{2 0 2}
-4 z^{3}_{0 1 1} \wedge z^{3}_{2 1 0} \wedge z^{4}_{0 1 2}
-2 z^{3}_{0 1 1} \wedge z^{3}_{2 1 0} \wedge z^{4}_{1 0 1}
-2 z^{3}_{0 1 1} \wedge z^{3}_{3 0 0} \wedge z^{4}_{0 1 1}
+4 z^{3}_{0 1 2} \wedge z^{3}_{0 2 0} \wedge z^{4}_{1 1 2} \\&
+3 z^{3}_{0 1 2} \wedge z^{3}_{0 2 1} \wedge z^{4}_{0 2 2}
-2 z^{3}_{0 1 2} \wedge z^{3}_{0 3 0} \wedge z^{4}_{0 1 3}
-4 z^{3}_{0 1 1} \wedge z^{3}_{1 2 0} \wedge z^{4}_{1 0 2}
+4 z^{3}_{0 1 1} \wedge z^{3}_{2 0 1} \wedge z^{4}_{0 2 1}
-2 z^{3}_{0 1 1} \wedge z^{3}_{2 0 1} \wedge z^{4}_{1 1 0} \\&
-z^{3}_{0 1 0} \wedge z^{3}_{2 0 1} \wedge z^{4}_{2 0 0}
-4 z^{3}_{0 1 0} \wedge z^{3}_{2 1 0} \wedge z^{4}_{1 0 2}
+4 z^{3}_{0 1 1} \wedge z^{3}_{0 2 0} \wedge z^{4}_{2 0 2}
+2 z^{3}_{0 1 1} \wedge z^{3}_{0 2 1} \wedge z^{4}_{1 1 2}
-2 z^{3}_{0 1 1} \wedge z^{3}_{0 3 0} \wedge z^{4}_{1 0 3} \\&
-2 z^{3}_{0 1 0} \wedge z^{3}_{3 0 0} \wedge z^{4}_{1 0 1}
+2 z^{3}_{0 1 1} \wedge z^{3}_{0 1 2} \wedge z^{4}_{1 2 1}
-2 z^{3}_{0 1 1} \wedge z^{3}_{1 0 0} \wedge z^{4}_{2 1 1}
-4 z^{3}_{0 1 1} \wedge z^{3}_{1 0 1} \wedge z^{4}_{1 2 1}
+4 z^{3}_{0 1 1} \wedge z^{3}_{1 0 1} \wedge z^{4}_{2 1 0} \\&
-2 z^{3}_{0 1 1} \wedge z^{3}_{1 0 2} \wedge z^{4}_{0 3 1}
+4 z^{3}_{0 1 1} \wedge z^{3}_{1 0 2} \wedge z^{4}_{1 2 0}
+4 z^{3}_{0 1 1} \wedge z^{3}_{1 1 0} \wedge z^{4}_{1 1 2}
+4 z^{3}_{0 1 1} \wedge z^{3}_{1 1 0} \wedge z^{4}_{2 0 1}
+4 z^{3}_{0 1 1} \wedge z^{3}_{1 1 1} \wedge z^{4}_{0 2 2} \\&
-2 z^{3}_{0 1 1} \wedge z^{3}_{1 2 0} \wedge z^{4}_{0 1 3}
+4 z^{3}_{0 1 0} \wedge z^{3}_{1 1 1} \wedge z^{4}_{1 1 2}
-2 z^{3}_{0 1 0} \wedge z^{3}_{1 1 1} \wedge z^{4}_{2 0 1}
-2 z^{3}_{0 1 0} \wedge z^{3}_{1 0 0} \wedge z^{4}_{3 0 1}
-4 z^{3}_{0 1 0} \wedge z^{3}_{1 0 1} \wedge z^{4}_{2 1 1} \\&
+2 z^{3}_{0 1 0} \wedge z^{3}_{1 0 1} \wedge z^{4}_{3 0 0}
-2 z^{3}_{0 1 0} \wedge z^{3}_{1 2 0} \wedge z^{4}_{1 0 3}
+4 z^{3}_{0 1 0} \wedge z^{3}_{2 0 0} \wedge z^{4}_{2 0 1}
+4 z^{3}_{0 1 0} \wedge z^{3}_{2 0 1} \wedge z^{4}_{1 1 1}
-2 z^{3}_{0 0 2} \wedge z^{3}_{2 0 1} \wedge z^{4}_{0 3 0} \\&
+2 z^{3}_{0 0 2} \wedge z^{3}_{2 1 0} \wedge z^{4}_{0 2 1}
+2 z^{3}_{0 0 2} \wedge z^{3}_{2 1 0} \wedge z^{4}_{1 1 0}
+z^{3}_{0 0 2} \wedge z^{3}_{3 0 0} \wedge z^{4}_{0 2 0}
+z^{3}_{0 0 3} \wedge z^{3}_{0 1 0} \wedge z^{4}_{2 2 0}
+2 z^{3}_{0 0 3} \wedge z^{3}_{0 1 1} \wedge z^{4}_{1 3 0} \\&
+z^{3}_{0 0 3} \wedge z^{3}_{0 1 2} \wedge z^{4}_{0 4 0}
-2 z^{3}_{0 0 3} \wedge z^{3}_{0 2 0} \wedge z^{4}_{1 2 1}
-2 z^{3}_{0 0 3} \wedge z^{3}_{0 2 1} \wedge z^{4}_{0 3 1}
+z^{3}_{0 0 3} \wedge z^{3}_{0 3 0} \wedge z^{4}_{0 2 2}
-2 z^{3}_{0 0 3} \wedge z^{3}_{1 1 0} \wedge z^{4}_{1 2 0} \\&
-2 z^{3}_{0 0 3} \wedge z^{3}_{1 1 1} \wedge z^{4}_{0 3 0}
+2 z^{3}_{0 0 3} \wedge z^{3}_{1 2 0} \wedge z^{4}_{0 2 1}
+z^{3}_{0 0 3} \wedge z^{3}_{2 1 0} \wedge z^{4}_{0 2 0}
+2 z^{3}_{0 1 0} \wedge z^{3}_{0 1 1} \wedge z^{4}_{3 0 1}
-2 z^{3}_{1 0 1} \wedge z^{3}_{1 0 2} \wedge z^{4}_{0 3 0} \\&
+z^{3}_{0 0 0} \wedge z^{3}_{1 0 0} \wedge z^{4}_{4 0 0} 
+2 z^{3}_{0 0 0} \wedge z^{3}_{1 0 1} \wedge z^{4}_{3 1 0}
+z^{3}_{0 0 0} \wedge z^{3}_{1 0 2} \wedge z^{4}_{2 2 0}
-2 z^{3}_{0 0 0} \wedge z^{3}_{1 1 0} \wedge z^{4}_{3 0 1}
-2 z^{3}_{0 0 0} \wedge z^{3}_{1 1 1} \wedge z^{4}_{2 1 1} \\&
+z^{3}_{0 0 0} \wedge z^{3}_{1 2 0} \wedge z^{4}_{2 0 2}
-2 z^{3}_{0 0 0} \wedge z^{3}_{2 0 0} \wedge z^{4}_{3 0 0}
-2 z^{3}_{0 0 0} \wedge z^{3}_{2 0 1} \wedge z^{4}_{2 1 0}
+2 z^{3}_{0 0 0} \wedge z^{3}_{2 1 0} \wedge z^{4}_{2 0 1}
+z^{3}_{0 0 0} \wedge z^{3}_{3 0 0} \wedge z^{4}_{2 0 0} \\&
+z^{3}_{0 0 1} \wedge z^{3}_{0 1 0} \wedge z^{4}_{4 0 0}
+2 z^{3}_{0 0 1} \wedge z^{3}_{0 1 1} \wedge z^{4}_{3 1 0}
+z^{3}_{0 0 1} \wedge z^{3}_{0 1 2} \wedge z^{4}_{2 2 0}
\end{align*}

%
We also prepare a basis of 
$\ds \CGF{4}{4}$: 
\begin{align*}  & B_1 =  \\
&
2 z^{3}_{100} \wedge z^{3}_{101} \wedge z^{3}_{110} \wedge z^{3}_{200}
+z^{3}_{020} \wedge z^{3}_{100} \wedge z^{3}_{102} \wedge z^{3}_{200}
+2 z^{3}_{020} \wedge z^{3}_{101} \wedge z^{3}_{102} \wedge z^{3}_{110}
+\frac{1}{3} z^{3}_{000} \wedge z^{3}_{001} \wedge z^{3}_{210} \wedge z^{3}_{300}\\& 
+z^{3}_{002} \wedge z^{3}_{100} \wedge z^{3}_{120} \wedge z^{3}_{200}
-2 z^{3}_{002} \wedge z^{3}_{101} \wedge z^{3}_{110} \wedge z^{3}_{120}
+\frac{1}{3} z^{3}_{003} \wedge z^{3}_{010} \wedge z^{3}_{030} \wedge z^{3}_{201}
+\frac{2}{3} z^{3}_{003} \wedge z^{3}_{011} \wedge z^{3}_{030} \wedge z^{3}_{111}\\& 
-\frac{1}{3} z^{3}_{003} \wedge z^{3}_{012} \wedge z^{3}_{021} \wedge z^{3}_{030}
-2 z^{3}_{002} \wedge z^{3}_{011} \wedge z^{3}_{111} \wedge z^{3}_{120}
-z^{3}_{002} \wedge z^{3}_{012} \wedge z^{3}_{021} \wedge z^{3}_{120}
+z^{3}_{002} \wedge z^{3}_{020} \wedge z^{3}_{102} \wedge z^{3}_{120}\\&
+z^{3}_{001} \wedge z^{3}_{020} \wedge z^{3}_{102} \wedge z^{3}_{210}
-z^{3}_{001} \wedge z^{3}_{100} \wedge z^{3}_{200} \wedge z^{3}_{210}
-2 z^{3}_{001} \wedge z^{3}_{101} \wedge z^{3}_{110} \wedge z^{3}_{210}
+z^{3}_{002} \wedge z^{3}_{010} \wedge z^{3}_{120} \wedge z^{3}_{201}\\&
+z^{3}_{001} \wedge z^{3}_{002} \wedge z^{3}_{120} \wedge z^{3}_{210}
+\frac{1}{3} z^{3}_{001} \wedge z^{3}_{003} \wedge z^{3}_{030} \wedge z^{3}_{210}
-z^{3}_{001} \wedge z^{3}_{010} \wedge z^{3}_{201} \wedge z^{3}_{210}
-2 z^{3}_{001} \wedge z^{3}_{011} \wedge z^{3}_{111} \wedge
z^{3}_{210}\\ \allowbreak &
-z^{3}_{001} \wedge z^{3}_{012} \wedge z^{3}_{021} \wedge z^{3}_{210}
+\frac{1}{3} z^{3}_{002} \wedge z^{3}_{003} \wedge z^{3}_{030} \wedge z^{3}_{120}
-2 z^{3}_{011} \wedge z^{3}_{020} \wedge z^{3}_{102} \wedge z^{3}_{111}
-2 z^{3}_{011} \wedge z^{3}_{100} \wedge z^{3}_{111} \wedge z^{3}_{200}\\&
+4 z^{3}_{011} \wedge z^{3}_{101} \wedge z^{3}_{110} \wedge z^{3}_{111}
+z^{3}_{012} \wedge z^{3}_{020} \wedge z^{3}_{021} \wedge z^{3}_{102}
+z^{3}_{012} \wedge z^{3}_{021} \wedge z^{3}_{100} \wedge z^{3}_{200}
+2 z^{3}_{012} \wedge z^{3}_{021} \wedge z^{3}_{101} \wedge z^{3}_{110}\\&
-\frac{1}{3} z^{3}_{003} \wedge z^{3}_{020} \wedge z^{3}_{030} \wedge z^{3}_{102}
+z^{3}_{010} \wedge z^{3}_{100} \wedge z^{3}_{200} \wedge z^{3}_{201}
+2 z^{3}_{010} \wedge z^{3}_{101} \wedge z^{3}_{110} \wedge z^{3}_{201}
+2 z^{3}_{011} \wedge z^{3}_{012} \wedge z^{3}_{021} \wedge z^{3}_{111}\\&
-\frac{2}{3} z^{3}_{003} \wedge z^{3}_{030} \wedge z^{3}_{101} \wedge z^{3}_{110}
+2 z^{3}_{010} \wedge z^{3}_{011} \wedge z^{3}_{111} \wedge z^{3}_{201}
+z^{3}_{010} \wedge z^{3}_{012} \wedge z^{3}_{021} \wedge z^{3}_{201}
-z^{3}_{010} \wedge z^{3}_{020} \wedge z^{3}_{102} \wedge z^{3}_{201}\\&
-\frac{1}{3} z^{3}_{003} \wedge z^{3}_{030} \wedge z^{3}_{100} \wedge z^{3}_{200}
+\frac{1}{3} z^{3}_{000} \wedge z^{3}_{002} \wedge z^{3}_{120} \wedge z^{3}_{300}
+\frac{1}{9} z^{3}_{000} \wedge z^{3}_{003} \wedge z^{3}_{030} \wedge z^{3}_{300}
-\frac{1}{3} z^{3}_{000} \wedge z^{3}_{010} \wedge z^{3}_{201} \wedge z^{3}_{300}\\&
-\frac{2}{3} z^{3}_{000} \wedge z^{3}_{011} \wedge z^{3}_{111} \wedge z^{3}_{300}
-\frac{1}{3} z^{3}_{000} \wedge z^{3}_{012} \wedge z^{3}_{021} \wedge z^{3}_{300}
+\frac{1}{3} z^{3}_{000} \wedge z^{3}_{020} \wedge z^{3}_{102} \wedge z^{3}_{300}
-\frac{1}{3} z^{3}_{000} \wedge z^{3}_{100} \wedge z^{3}_{200} \wedge z^{3}_{300}\\&
-\frac{2}{3} z^{3}_{000} \wedge z^{3}_{101} \wedge z^{3}_{110} \wedge z^{3}_{300}
\end{align*}

\begin{align*}  & B_2 =  \\&
\frac{1}{2} z^{3}_{002} \wedge z^{3}_{030} \wedge z^{3}_{100} \wedge z^{3}_{201}
-z^{3}_{002} \wedge z^{3}_{021} \wedge z^{3}_{110} \wedge z^{3}_{200}
-\frac{1}{2} z^{3}_{002} \wedge z^{3}_{021} \wedge z^{3}_{100} \wedge z^{3}_{210}
-z^{3}_{100} \wedge z^{3}_{101} \wedge z^{3}_{110} \wedge z^{3}_{200} \\&
+2 z^{3}_{021} \wedge z^{3}_{100} \wedge z^{3}_{101} \wedge z^{3}_{111}
-z^{3}_{021} \wedge z^{3}_{100} \wedge z^{3}_{102} \wedge z^{3}_{110}
-z^{3}_{030} \wedge z^{3}_{100} \wedge z^{3}_{101} \wedge z^{3}_{102}
-z^{3}_{020} \wedge z^{3}_{100} \wedge z^{3}_{102} \wedge z^{3}_{111} \\&
-\frac{1}{2} z^{3}_{020} \wedge z^{3}_{100} \wedge z^{3}_{102} \wedge z^{3}_{200}
+z^{3}_{020} \wedge z^{3}_{101} \wedge z^{3}_{102} \wedge z^{3}_{110}
-\frac{1}{2} z^{3}_{001} \wedge z^{3}_{030} \wedge z^{3}_{102} \wedge z^{3}_{200}
-\frac{1}{6} z^{3}_{000} \wedge z^{3}_{001} \wedge z^{3}_{210} \wedge z^{3}_{300} \\&
+z^{3}_{002} \wedge z^{3}_{030} \wedge z^{3}_{101} \wedge z^{3}_{200}
-z^{3}_{012} \wedge z^{3}_{100} \wedge z^{3}_{101} \wedge z^{3}_{120}
-2 z^{3}_{012} \wedge z^{3}_{100}z^{3}_{110} \wedge z^{3}_{111}
-z^{3}_{002} \wedge z^{3}_{100} \wedge z^{3}_{111} \wedge z^{3}_{120} \\&
-\frac{1}{2} z^{3}_{002} \wedge z^{3}_{100} \wedge z^{3}_{120} \wedge z^{3}_{200}
-z^{3}_{002} \wedge z^{3}_{101} \wedge z^{3}_{110} \wedge z^{3}_{120}
-\frac{1}{6} z^{3}_{003} \wedge z^{3}_{010} \wedge z^{3}_{030} \wedge z^{3}_{201}
-z^{3}_{003} \wedge z^{3}_{010} \wedge z^{3}_{110} \wedge z^{3}_{210} \\&
-\frac{1}{2} z^{3}_{003} \wedge z^{3}_{010} \wedge z^{3}_{120} \wedge z^{3}_{200}
+\frac{1}{2} z^{3}_{003} \wedge z^{3}_{010} \wedge z^{3}_{020} \wedge z^{3}_{300}
-\frac{1}{3} z^{3}_{003} \wedge z^{3}_{011} \wedge z^{3}_{030} \wedge z^{3}_{111}
+\frac{1}{6} z^{3}_{003} \wedge z^{3}_{012} \wedge z^{3}_{021} \wedge z^{3}_{030} \\&
+2 z^{3}_{002} \wedge z^{3}_{011} \wedge z^{3}_{110} \wedge z^{3}_{210}
+z^{3}_{002} \wedge z^{3}_{011} \wedge z^{3}_{111} \wedge z^{3}_{120}
+z^{3}_{002} \wedge z^{3}_{011} \wedge z^{3}_{120} \wedge z^{3}_{200}
+\frac{1}{2} z^{3}_{002} \wedge z^{3}_{012} \wedge z^{3}_{021} \wedge z^{3}_{120} \\&
-z^{3}_{002} \wedge z^{3}_{020} \wedge z^{3}_{110} \wedge z^{3}_{201}
-z^{3}_{002} \wedge z^{3}_{020} \wedge z^{3}_{111} \wedge z^{3}_{200}
-z^{3}_{002} \wedge z^{3}_{020} \wedge z^{3}_{101} \wedge z^{3}_{210}
-\frac{1}{2} z^{3}_{002} \wedge z^{3}_{020} \wedge z^{3}_{102} \wedge z^{3}_{120} \\&
+\frac{1}{2} z^{3}_{001} \wedge z^{3}_{100} \wedge z^{3}_{200} \wedge z^{3}_{210}
+z^{3}_{001} \wedge z^{3}_{101} \wedge z^{3}_{110} \wedge z^{3}_{210}
+2 z^{3}_{001} \wedge z^{3}_{101} \wedge z^{3}_{111} \wedge z^{3}_{120}
-\frac{1}{2} z^{3}_{002} \wedge z^{3}_{010} \wedge z^{3}_{021} \wedge z^{3}_{300} \\&
+z^{3}_{002} \wedge z^{3}_{010} \wedge z^{3}_{111} \wedge z^{3}_{210}
-z^{3}_{002} \wedge z^{3}_{011} \wedge z^{3}_{020} \wedge z^{3}_{300}
+z^{3}_{001} \wedge z^{3}_{102} \wedge z^{3}_{110} \wedge z^{3}_{120}
-z^{3}_{000} \wedge z^{3}_{102} \wedge z^{3}_{111} \wedge z^{3}_{120} \\&
-\frac{1}{2} z^{3}_{001} \wedge z^{3}_{002} \wedge z^{3}_{030} \wedge z^{3}_{300}
+z^{3}_{001} \wedge z^{3}_{002} \wedge z^{3}_{120} \wedge z^{3}_{210}
-\frac{1}{6} z^{3}_{001} \wedge z^{3}_{003} \wedge z^{3}_{030} \wedge z^{3}_{210}
+\frac{1}{2} z^{3}_{001} \wedge z^{3}_{010} \wedge z^{3}_{201} \wedge z^{3}_{210} \\&
+z^{3}_{001} \wedge z^{3}_{011}z^{3}_{021} \wedge z^{3}_{300}
-z^{3}_{001} \wedge z^{3}_{011} \wedge z^{3}_{111} \wedge z^{3}_{210}
-z^{3}_{001} \wedge z^{3}_{011} \wedge z^{3}_{120} \wedge z^{3}_{201}
+\frac{1}{2} z^{3}_{001} \wedge z^{3}_{012} \wedge z^{3}_{020} \wedge z^{3}_{300} \\&
+\frac{1}{2} z^{3}_{001} \wedge z^{3}_{012} \wedge z^{3}_{021} \wedge z^{3}_{210}
-z^{3}_{001} \wedge z^{3}_{012} \wedge z^{3}_{110} \wedge z^{3}_{210}
-\frac{1}{2} z^{3}_{001} \wedge z^{3}_{012} \wedge z^{3}_{120} \wedge z^{3}_{200}
+z^{3}_{001} \wedge z^{3}_{020} \wedge z^{3}_{111} \wedge z^{3}_{201} \\&
+z^{3}_{001} \wedge z^{3}_{021} \wedge z^{3}_{101} \wedge z^{3}_{210}
+z^{3}_{001} \wedge z^{3}_{021} \wedge z^{3}_{110} \wedge z^{3}_{201}
+z^{3}_{001} \wedge z^{3}_{021} \wedge z^{3}_{111} \wedge z^{3}_{200}
-z^{3}_{001} \wedge z^{3}_{030} \wedge z^{3}_{101} \wedge z^{3}_{201} \\&
-\frac{1}{6} z^{3}_{002} \wedge z^{3}_{003} \wedge z^{3}_{030} \wedge z^{3}_{120}
+z^{3}_{011} \wedge z^{3}_{020} \wedge z^{3}_{102} \wedge z^{3}_{111}
+z^{3}_{011} \wedge z^{3}_{020} \wedge z^{3}_{102} \wedge z^{3}_{200}
-z^{3}_{011} \wedge z^{3}_{021} \wedge z^{3}_{100} \wedge z^{3}_{201} \\&
-2 z^{3}_{011} \wedge z^{3}_{021} \wedge z^{3}_{101} \wedge z^{3}_{200}
+z^{3}_{011} \wedge z^{3}_{100} \wedge z^{3}_{102} \wedge z^{3}_{120}
+z^{3}_{011} \wedge z^{3}_{100} \wedge z^{3}_{111} \wedge z^{3}_{200}
+2 z^{3}_{011} \wedge z^{3}_{101} \wedge z^{3}_{110} \wedge z^{3}_{111} \\&
-\frac{1}{2} z^{3}_{012} \wedge z^{3}_{020} \wedge z^{3}_{021} \wedge z^{3}_{102}
-\frac{1}{2} z^{3}_{012} \wedge z^{3}_{020} \wedge z^{3}_{100} \wedge z^{3}_{201}
-z^{3}_{012} \wedge z^{3}_{020} \wedge z^{3}_{101} \wedge z^{3}_{200}
+z^{3}_{012} \wedge z^{3}_{021} \wedge z^{3}_{100} \wedge z^{3}_{200} \\&
-z^{3}_{012} \wedge z^{3}_{021} \wedge z^{3}_{101} \wedge z^{3}_{110}
+\frac{1}{6} z^{3}_{003} \wedge z^{3}_{020} \wedge z^{3}_{030} \wedge z^{3}_{102}
+\frac{1}{2} z^{3}_{003} \wedge z^{3}_{020} \wedge z^{3}_{100} \wedge z^{3}_{210}
-\frac{1}{2} z^{3}_{010} \wedge z^{3}_{100} \wedge z^{3}_{200} \wedge z^{3}_{201} \\&
-z^{3}_{010} \wedge z^{3}_{101} \wedge z^{3}_{102} \wedge z^{3}_{120}
-z^{3}_{010} \wedge z^{3}_{101} \wedge z^{3}_{110} \wedge z^{3}_{201}
-2 z^{3}_{010} \wedge z^{3}_{102} \wedge z^{3}_{110} \wedge z^{3}_{111}
-z^{3}_{011} \wedge z^{3}_{012} \wedge z^{3}_{021} \wedge z^{3}_{111} \\&
+z^{3}_{011} \wedge z^{3}_{012} \wedge z^{3}_{100} \wedge z^{3}_{210}
+2 z^{3}_{011} \wedge z^{3}_{012} \wedge z^{3}_{110} \wedge z^{3}_{200}
+2 z^{3}_{011} \wedge z^{3}_{020} \wedge z^{3}_{101} \wedge z^{3}_{201}
+\frac{1}{3} z^{3}_{003} \wedge z^{3}_{030} \wedge z^{3}_{101} \wedge z^{3}_{110} \\&
+z^{3}_{003} \wedge z^{3}_{100} \wedge z^{3}_{110} \wedge z^{3}_{120}
-z^{3}_{010} \wedge z^{3}_{011} \wedge z^{3}_{012} \wedge z^{3}_{300}
-z^{3}_{010} \wedge z^{3}_{012} \wedge z^{3}_{101} \wedge z^{3}_{210}
-z^{3}_{010} \wedge z^{3}_{012} \wedge z^{3}_{110} \wedge z^{3}_{201} \\&
-z^{3}_{010} \wedge z^{3}_{012} \wedge z^{3}_{111} \wedge z^{3}_{200}
+z^{3}_{010} \wedge z^{3}_{011} \wedge z^{3}_{102} \wedge z^{3}_{210}
+z^{3}_{010} \wedge z^{3}_{011} \wedge z^{3}_{111} \wedge z^{3}_{201}
-\frac{1}{2} z^{3}_{010} \wedge z^{3}_{012} \wedge z^{3}_{021} \wedge z^{3}_{201} \\&
-z^{3}_{010} \wedge z^{3}_{020} \wedge z^{3}_{102} \wedge z^{3}_{201}
+z^{3}_{010} \wedge z^{3}_{021} \wedge z^{3}_{101} \wedge z^{3}_{201}
+\frac{1}{2} z^{3}_{010} \wedge z^{3}_{021} \wedge z^{3}_{102} \wedge z^{3}_{200}
+z^{3}_{003} \wedge z^{3}_{020} \wedge z^{3}_{110} \wedge z^{3}_{200} \\&
-\frac{1}{3} z^{3}_{003} \wedge z^{3}_{030} \wedge z^{3}_{100} \wedge z^{3}_{200}
-\frac{1}{6} z^{3}_{000} \wedge z^{3}_{002} \wedge z^{3}_{120} \wedge z^{3}_{300}
+\frac{1}{9} z^{3}_{000} \wedge z^{3}_{003} \wedge z^{3}_{030} \wedge z^{3}_{300}
-\frac{1}{2} z^{3}_{000} \wedge z^{3}_{003} \wedge z^{3}_{120} \wedge z^{3}_{210} \\&
+\frac{1}{6} z^{3}_{000} \wedge z^{3}_{010} \wedge z^{3}_{201} \wedge z^{3}_{300}
+\frac{1}{3} z^{3}_{000} \wedge z^{3}_{011} \wedge z^{3}_{111} \wedge z^{3}_{300}
-\frac{1}{3} z^{3}_{000} \wedge z^{3}_{012} \wedge z^{3}_{021} \wedge z^{3}_{300}
+z^{3}_{000} \wedge z^{3}_{012} \wedge z^{3}_{111} \wedge z^{3}_{210} \\&
+\frac{1}{2} z^{3}_{000} \wedge z^{3}_{012} \wedge z^{3}_{120} \wedge z^{3}_{201}
-\frac{1}{6} z^{3}_{000} \wedge z^{3}_{020} \wedge z^{3}_{102} \wedge z^{3}_{300}
-\frac{1}{2} z^{3}_{000} \wedge z^{3}_{021} \wedge z^{3}_{102} \wedge z^{3}_{210}
-z^{3}_{000} \wedge z^{3}_{021} \wedge z^{3}_{111} \wedge z^{3}_{201} \\&
+\frac{1}{2} z^{3}_{000} \wedge z^{3}_{030} \wedge z^{3}_{102} \wedge z^{3}_{201}
+\frac{1}{6} z^{3}_{000} \wedge z^{3}_{100} \wedge z^{3}_{200} \wedge z^{3}_{300}
+\frac{1}{3} z^{3}_{000} \wedge z^{3}_{101} \wedge z^{3}_{110} \wedge z^{3}_{300}
\end{align*} 

\begin{align*} & B_3 =  \\&
\frac{1}{2} z^{3}_{002} \wedge z^{3}_{030}z^{3}_{100} \wedge z^{3}_{201}
-z^{3}_{002} \wedge z^{3}_{021} \wedge z^{3}_{110} \wedge z^{3}_{200}
-\frac{1}{2} z^{3}_{002} \wedge z^{3}_{021} \wedge z^{3}_{100} \wedge z^{3}_{210}
+z^{3}_{002} \wedge z^{3}_{021} \wedge z^{3}_{101} \wedge z^{3}_{120} \\&
-4 z^{3}_{100} \wedge z^{3}_{101} \wedge z^{3}_{110} \wedge z^{3}_{200}
-z^{3}_{002} \wedge z^{3}_{021} \wedge z^{3}_{030} \wedge z^{3}_{102}
+2 z^{3}_{021} \wedge z^{3}_{100} \wedge z^{3}_{101} \wedge z^{3}_{111}
-z^{3}_{021} \wedge z^{3}_{100} \wedge z^{3}_{102} \wedge z^{3}_{110} \\&
-z^{3}_{030} \wedge z^{3}_{100} \wedge z^{3}_{101} \wedge z^{3}_{102}
-z^{3}_{020} \wedge z^{3}_{100} \wedge z^{3}_{101} \wedge z^{3}_{201}
-2 z^{3}_{020} \wedge z^{3}_{100} \wedge z^{3}_{102} \wedge z^{3}_{111}
-\frac{3}{2} z^{3}_{020} \wedge z^{3}_{100} \wedge z^{3}_{102} \wedge z^{3}_{200} \\&
-\frac{1}{2} z^{3}_{001} \wedge z^{3}_{030} \wedge z^{3}_{102} \wedge z^{3}_{200}
-\frac{1}{2} z^{3}_{000} \wedge z^{3}_{001} \wedge z^{3}_{210} \wedge z^{3}_{300}
-z^{3}_{002} \wedge z^{3}_{030} \wedge z^{3}_{101} \wedge z^{3}_{111}
+z^{3}_{002} \wedge z^{3}_{030} \wedge z^{3}_{101} \wedge z^{3}_{200} \\&
-z^{3}_{012} \wedge z^{3}_{100} \wedge z^{3}_{101} \wedge z^{3}_{120}
-2 z^{3}_{012} \wedge z^{3}_{100} \wedge z^{3}_{110} \wedge z^{3}_{111}
-z^{3}_{020} \wedge z^{3}_{021} \wedge z^{3}_{101} \wedge z^{3}_{102}
-2 z^{3}_{002} \wedge z^{3}_{100} \wedge z^{3}_{111} \wedge z^{3}_{120} \\&
-\frac{3}{2} z^{3}_{002} \wedge z^{3}_{100} \wedge z^{3}_{120} \wedge z^{3}_{200}
-\frac{1}{2} z^{3}_{003} \wedge z^{3}_{010} \wedge z^{3}_{030} \wedge z^{3}_{201}
-2 z^{3}_{003} \wedge z^{3}_{010} \wedge z^{3}_{110} \wedge z^{3}_{210}
-z^{3}_{003} \wedge z^{3}_{010} \wedge z^{3}_{111} \wedge z^{3}_{120} \\&
-\frac{1}{2} z^{3}_{003} \wedge z^{3}_{010} \wedge z^{3}_{120} \wedge z^{3}_{200}
-z^{3}_{003} \wedge z^{3}_{011} \wedge z^{3}_{020} \wedge z^{3}_{210}
+\frac{1}{2} z^{3}_{003} \wedge z^{3}_{010} \wedge z^{3}_{020} \wedge z^{3}_{300}
-z^{3}_{003} \wedge z^{3}_{011} \wedge z^{3}_{021} \wedge z^{3}_{120} \\&
-z^{3}_{003} \wedge z^{3}_{011} \wedge z^{3}_{110} \wedge z^{3}_{120}
-z^{3}_{003} \wedge z^{3}_{012} \wedge z^{3}_{020} \wedge z^{3}_{120}
-\frac{1}{2} z^{3}_{003} \wedge z^{3}_{012} \wedge z^{3}_{021} \wedge z^{3}_{030}
+z^{3}_{002} \wedge z^{3}_{011} \wedge z^{3}_{030} \wedge z^{3}_{201} \\&
+z^{3}_{002} \wedge z^{3}_{011} \wedge z^{3}_{111} \wedge z^{3}_{120}
+2 z^{3}_{002} \wedge z^{3}_{011} \wedge z^{3}_{120} \wedge z^{3}_{200}
-z^{3}_{002} \wedge z^{3}_{012} \wedge z^{3}_{020} \wedge z^{3}_{210}
+\frac{1}{2} z^{3}_{002} \wedge z^{3}_{012} \wedge z^{3}_{021} \wedge z^{3}_{120} \\&
+z^{3}_{002} \wedge z^{3}_{012} \wedge z^{3}_{030} \wedge z^{3}_{111}
-2 z^{3}_{002} \wedge z^{3}_{020} \wedge z^{3}_{111} \wedge z^{3}_{200}
-z^{3}_{002} \wedge z^{3}_{012} \wedge z^{3}_{110} \wedge z^{3}_{120}
-z^{3}_{002} \wedge z^{3}_{020} \wedge z^{3}_{021} \wedge z^{3}_{201} \\&
-\frac{1}{2} z^{3}_{002} \wedge z^{3}_{020} \wedge z^{3}_{102} \wedge z^{3}_{120}
-z^{3}_{002} \wedge z^{3}_{100} \wedge z^{3}_{110} \wedge z^{3}_{210}
+z^{3}_{001} \wedge z^{3}_{101} \wedge z^{3}_{120} \wedge z^{3}_{200}
-z^{3}_{001} \wedge z^{3}_{110} \wedge z^{3}_{111} \wedge z^{3}_{200} \\&
-z^{3}_{001} \wedge z^{3}_{100} \wedge z^{3}_{110} \wedge z^{3}_{300}
-z^{3}_{001} \wedge z^{3}_{100} \wedge z^{3}_{111} \wedge z^{3}_{210}
+z^{3}_{001} \wedge z^{3}_{100} \wedge z^{3}_{120} \wedge z^{3}_{201}
+\frac{1}{2} z^{3}_{001} \wedge z^{3}_{100} \wedge z^{3}_{200} \wedge z^{3}_{210} \\&
+z^{3}_{001} \wedge z^{3}_{101} \wedge z^{3}_{110} \wedge z^{3}_{210}
-\frac{1}{2} z^{3}_{002} \wedge z^{3}_{010} \wedge z^{3}_{021} \wedge z^{3}_{300}
-z^{3}_{002} \wedge z^{3}_{010} \wedge z^{3}_{110} \wedge z^{3}_{300}
-z^{3}_{002} \wedge z^{3}_{010} \wedge z^{3}_{200} \wedge z^{3}_{210} \\&
-2 z^{3}_{002} \wedge z^{3}_{011} \wedge z^{3}_{020} \wedge z^{3}_{300}
-z^{3}_{002} \wedge z^{3}_{011} \wedge z^{3}_{021} \wedge z^{3}_{210}
-z^{3}_{000} \wedge z^{3}_{101} \wedge z^{3}_{111} \wedge z^{3}_{210}
+z^{3}_{000} \wedge z^{3}_{101} \wedge z^{3}_{120} \wedge z^{3}_{201} \\&
-z^{3}_{000} \wedge z^{3}_{101} \wedge z^{3}_{200} \wedge z^{3}_{210}
-z^{3}_{000} \wedge z^{3}_{102} \wedge z^{3}_{110} \wedge z^{3}_{210}
-2 z^{3}_{000} \wedge z^{3}_{102} \wedge z^{3}_{111} \wedge z^{3}_{120}
-z^{3}_{000} \wedge z^{3}_{110} \wedge z^{3}_{111} \wedge z^{3}_{201} \\&
+z^{3}_{000} \wedge z^{3}_{110} \wedge z^{3}_{200} \wedge z^{3}_{201}
-\frac{1}{2} z^{3}_{001} \wedge z^{3}_{002} \wedge z^{3}_{030} \wedge z^{3}_{300}
-\frac{1}{2} z^{3}_{001} \wedge z^{3}_{003} \wedge z^{3}_{030} \wedge z^{3}_{210}
-z^{3}_{001} \wedge z^{3}_{010} \wedge z^{3}_{200} \wedge z^{3}_{300} \\&
+\frac{1}{2} z^{3}_{001} \wedge z^{3}_{010} \wedge z^{3}_{201} \wedge z^{3}_{210}
+z^{3}_{001} \wedge z^{3}_{011} \wedge z^{3}_{021} \wedge z^{3}_{300}
-z^{3}_{001} \wedge z^{3}_{011} \wedge z^{3}_{110} \wedge z^{3}_{300}
-z^{3}_{001} \wedge z^{3}_{011} \wedge z^{3}_{200} \wedge z^{3}_{210} \\&
+\frac{1}{2} z^{3}_{001} \wedge z^{3}_{012} \wedge z^{3}_{020} \wedge z^{3}_{300}
+\frac{3}{2} z^{3}_{001} \wedge z^{3}_{012} \wedge z^{3}_{021} \wedge z^{3}_{210}
-2 z^{3}_{001} \wedge z^{3}_{012} \wedge z^{3}_{110} \wedge z^{3}_{210}
-z^{3}_{001} \wedge z^{3}_{012} \wedge z^{3}_{111} \wedge z^{3}_{120} \\&
-\frac{1}{2} z^{3}_{001} \wedge z^{3}_{012} \wedge z^{3}_{120} \wedge z^{3}_{200}
+z^{3}_{001} \wedge z^{3}_{020} \wedge z^{3}_{101} \wedge z^{3}_{300}
+z^{3}_{001} \wedge z^{3}_{020} \wedge z^{3}_{200} \wedge z^{3}_{201}
+2 z^{3}_{001} \wedge z^{3}_{021} \wedge z^{3}_{101} \wedge z^{3}_{210} \\&
+z^{3}_{001} \wedge z^{3}_{021}z^{3}_{102} \wedge z^{3}_{120}
+2 z^{3}_{001} \wedge z^{3}_{021} \wedge z^{3}_{110} \wedge z^{3}_{201}
+z^{3}_{001} \wedge z^{3}_{021} \wedge z^{3}_{111} \wedge z^{3}_{200}
-2 z^{3}_{001} \wedge z^{3}_{030}z^{3}_{101} \wedge z^{3}_{201} \\&
-z^{3}_{001} \wedge z^{3}_{030} \wedge z^{3}_{102} \wedge z^{3}_{111}
-\frac{1}{2} z^{3}_{002} \wedge z^{3}_{003} \wedge z^{3}_{030} \wedge z^{3}_{120}
+z^{3}_{011} \wedge z^{3}_{020} \wedge z^{3}_{102} \wedge z^{3}_{111}
+2 z^{3}_{011} \wedge z^{3}_{020} \wedge z^{3}_{102} \wedge z^{3}_{200} \\&
-z^{3}_{011} \wedge z^{3}_{021} \wedge z^{3}_{100} \wedge z^{3}_{201}
-2 z^{3}_{011} \wedge z^{3}_{021} \wedge z^{3}_{101} \wedge z^{3}_{200}
+z^{3}_{011} \wedge z^{3}_{021} \wedge z^{3}_{102} \wedge z^{3}_{110}
+z^{3}_{011} \wedge z^{3}_{030} \wedge z^{3}_{101} \wedge z^{3}_{102} \\&
+z^{3}_{011} \wedge z^{3}_{100} \wedge z^{3}_{101} \wedge z^{3}_{210}
+2 z^{3}_{011} \wedge z^{3}_{100} \wedge z^{3}_{102} \wedge z^{3}_{120}
+z^{3}_{011} \wedge z^{3}_{100} \wedge z^{3}_{110} \wedge z^{3}_{201}
+3 z^{3}_{011} \wedge z^{3}_{100} \wedge z^{3}_{111} \wedge z^{3}_{200} \\&
-\frac{1}{2} z^{3}_{012} \wedge z^{3}_{020} \wedge z^{3}_{021} \wedge z^{3}_{102}
-\frac{1}{2} z^{3}_{012} \wedge z^{3}_{020} \wedge z^{3}_{100} \wedge z^{3}_{201}
-z^{3}_{012} \wedge z^{3}_{020} \wedge z^{3}_{101} \wedge z^{3}_{111}
-z^{3}_{012} \wedge z^{3}_{020} \wedge z^{3}_{101} \wedge z^{3}_{200} \\&
+z^{3}_{012} \wedge z^{3}_{020} \wedge z^{3}_{102} \wedge z^{3}_{110}
-3 z^{3}_{012} \wedge z^{3}_{021} \wedge z^{3}_{101} \wedge z^{3}_{110}
-z^{3}_{003} \wedge z^{3}_{020} \wedge z^{3}_{021} \wedge z^{3}_{111}
+\frac{1}{2} z^{3}_{003} \wedge z^{3}_{020} \wedge z^{3}_{030} \wedge z^{3}_{102} \\&
+\frac{1}{2} z^{3}_{003} \wedge z^{3}_{020} \wedge z^{3}_{100} \wedge z^{3}_{210}
+z^{3}_{003} \wedge z^{3}_{020} \wedge z^{3}_{110} \wedge z^{3}_{111}
-z^{3}_{010} \wedge z^{3}_{100} \wedge z^{3}_{111} \wedge z^{3}_{201}
-\frac{1}{2} z^{3}_{010} \wedge z^{3}_{100} \wedge z^{3}_{200} \wedge z^{3}_{201} \\&
-z^{3}_{010} \wedge z^{3}_{101} \wedge z^{3}_{110} \wedge z^{3}_{201}
-z^{3}_{010} \wedge z^{3}_{101} \wedge z^{3}_{111} \wedge z^{3}_{200}
-z^{3}_{010} \wedge z^{3}_{102} \wedge z^{3}_{110} \wedge z^{3}_{200}
-z^{3}_{011} \wedge z^{3}_{012} \wedge z^{3}_{020} \wedge z^{3}_{201}   \\&
-4 z^{3}_{011} \wedge z^{3}_{012} \wedge z^{3}_{021} \wedge z^{3}_{111}
+z^{3}_{011} \wedge z^{3}_{012} \wedge z^{3}_{030} \wedge z^{3}_{102}
+z^{3}_{011} \wedge z^{3}_{012} \wedge z^{3}_{100} \wedge z^{3}_{210}
-z^{3}_{011} \wedge z^{3}_{012} \wedge z^{3}_{101} \wedge z^{3}_{120} \\&
+2 z^{3}_{011} \wedge z^{3}_{012} \wedge z^{3}_{110} \wedge z^{3}_{200}
+z^{3}_{003} \wedge z^{3}_{030} \wedge z^{3}_{101} \wedge z^{3}_{110}
+z^{3}_{003} \wedge z^{3}_{100} \wedge z^{3}_{110} \wedge z^{3}_{120}
-z^{3}_{010} \wedge z^{3}_{011} \wedge z^{3}_{012} \wedge z^{3}_{300} \\&
-2 z^{3}_{010} \wedge z^{3}_{012} \wedge z^{3}_{101} \wedge z^{3}_{210}
-z^{3}_{010} \wedge z^{3}_{012} \wedge z^{3}_{102} \wedge z^{3}_{120}
-2 z^{3}_{010} \wedge z^{3}_{012} \wedge z^{3}_{110} \wedge z^{3}_{201}
-z^{3}_{010} \wedge z^{3}_{012} \wedge z^{3}_{111} \wedge z^{3}_{200} \\&
-z^{3}_{010} \wedge z^{3}_{011} \wedge z^{3}_{101} \wedge z^{3}_{300}
-z^{3}_{010} \wedge z^{3}_{011} \wedge z^{3}_{200} \wedge z^{3}_{201}
-\frac{3}{2} z^{3}_{010} \wedge z^{3}_{012} \wedge z^{3}_{021} \wedge z^{3}_{201}
+2 z^{3}_{010} \wedge z^{3}_{021} \wedge z^{3}_{101} \wedge z^{3}_{201} \\&
+z^{3}_{010} \wedge z^{3}_{021} \wedge z^{3}_{102} \wedge z^{3}_{111}
+\frac{1}{2} z^{3}_{010} \wedge z^{3}_{021} \wedge z^{3}_{102} \wedge z^{3}_{200}
+z^{3}_{010} \wedge z^{3}_{100} \wedge z^{3}_{101} \wedge z^{3}_{300}
+z^{3}_{010} \wedge z^{3}_{100} \wedge z^{3}_{102} \wedge z^{3}_{210} \\&
+z^{3}_{003} \wedge z^{3}_{020} \wedge z^{3}_{110} \wedge z^{3}_{200}
+z^{3}_{002} \wedge z^{3}_{021} \wedge z^{3}_{110} \wedge z^{3}_{111}
-\frac{1}{2} z^{3}_{000} \wedge z^{3}_{002} \wedge z^{3}_{120} \wedge z^{3}_{300}
-\frac{1}{2} z^{3}_{000} \wedge z^{3}_{003} \wedge z^{3}_{120} \wedge z^{3}_{210} \\&
+\frac{1}{2} z^{3}_{000} \wedge z^{3}_{010} \wedge z^{3}_{201} \wedge z^{3}_{300}
+z^{3}_{000} \wedge z^{3}_{011} \wedge z^{3}_{111} \wedge z^{3}_{300} 
+z^{3}_{000} \wedge z^{3}_{012} \wedge z^{3}_{111} \wedge z^{3}_{210}
+\frac{1}{2} z^{3}_{000} \wedge z^{3}_{012} \wedge z^{3}_{120} \wedge z^{3}_{201} \\&
-\frac{1}{2} z^{3}_{000} \wedge z^{3}_{020} \wedge z^{3}_{102} \wedge z^{3}_{300}
-\frac{1}{2} z^{3}_{000} \wedge z^{3}_{021} \wedge z^{3}_{102} \wedge z^{3}_{210}
-z^{3}_{000} \wedge z^{3}_{021} \wedge z^{3}_{111} \wedge z^{3}_{201}
+\frac{1}{2} z^{3}_{000} \wedge z^{3}_{030} \wedge z^{3}_{102} \wedge z^{3}_{201} \\&
-\frac{1}{2} z^{3}_{000} \wedge z^{3}_{100} \wedge z^{3}_{200} \wedge z^{3}_{300}
-z^{3}_{000} \wedge z^{3}_{100} \wedge z^{3}_{201} \wedge z^{3}_{210} 
\end{align*}

\def\thesection{Appendix \Alph{section}}
\section{Concrete bases for weight $=6$}\label{Open:wt:six} 

\begin{align*}  B & := \frac{1}{2} z^{3}_{011}\wedge z^{4}_{220}\wedge z^{5}_{013}
-\frac{1}{6} z^{3}_{012}\wedge z^{4}_{010}\wedge z^{5}_{311}
-\frac{1}{12} z^{3}_{012}\wedge z^{4}_{004}\wedge z^{5}_{050} 
+\frac{1}{2} z^{3}_{011}\wedge z^{4}_{101}\wedge z^{5}_{310}
-z^{3}_{011}\wedge z^{4}_{120}\wedge z^{5}_{202}\\& 
+\frac{1}{2} z^{3}_{011}\wedge z^{4}_{112}\wedge z^{5}_{032} 
-\frac{1}{2} z^{3}_{011}\wedge z^{4}_{121}\wedge z^{5}_{023}
+\frac{1}{3} z^{3}_{011}\wedge z^{4}_{013}\wedge z^{5}_{131}
+\frac{1}{2} z^{3}_{011}\wedge z^{4}_{110}\wedge z^{5}_{212} 
+\frac{1}{6} z^{3}_{011}\wedge z^{4}_{301}\wedge z^{5}_{110}\\& 
+\frac{1}{2} z^{3}_{011}\wedge z^{4}_{130}\wedge z^{5}_{103}
-\frac{1}{4} z^{3}_{200}\wedge z^{4}_{310}\wedge z^{5}_{001} 
-\frac{1}{2} z^{3}_{011}\wedge z^{4}_{300}\wedge z^{5}_{111}
+\frac{1}{12} z^{3}_{200}\wedge z^{4}_{103}\wedge z^{5}_{030}
-\frac{1}{12} z^{3}_{210}\wedge z^{4}_{004}\wedge z^{5}_{030} \\& 
-\frac{1}{2} z^{3}_{011}\wedge z^{4}_{120}\wedge z^{5}_{113}
+\frac{1}{2} z^{3}_{011}\wedge z^{4}_{020}\wedge z^{5}_{302}
+\frac{1}{2} z^{3}_{011}\wedge z^{4}_{202}\wedge z^{5}_{031} 
-\frac{1}{6} z^{3}_{011}\wedge z^{4}_{004}\wedge z^{5}_{140}
-\frac{1}{2} z^{3}_{011}\wedge z^{4}_{201}\wedge z^{5}_{210}\\& 
-\frac{1}{6} z^{3}_{011}\wedge z^{4}_{103}\wedge z^{5}_{041} 
-z^{3}_{011}\wedge z^{4}_{210}\wedge z^{5}_{112}
+\frac{1}{2} z^{3}_{011}\wedge z^{4}_{200}\wedge z^{5}_{211}
-\frac{1}{2} z^{3}_{011}\wedge z^{4}_{301}\wedge z^{5}_{021} 
-z^{3}_{010}\wedge z^{4}_{211}\wedge z^{5}_{112}\\& 
+\frac{1}{6} z^{3}_{011}\wedge z^{4}_{040}\wedge z^{5}_{104}
-\frac{1}{6} z^{3}_{011}\wedge z^{4}_{100}\wedge z^{5}_{311} 
-\frac{1}{2} z^{3}_{011}\wedge z^{4}_{102}\wedge z^{5}_{131}
-\frac{1}{2} z^{3}_{011}\wedge z^{4}_{121}\wedge z^{5}_{112}
-\frac{1}{12} z^{3}_{012}\wedge z^{4}_{001}\wedge z^{5}_{320} \\& 
-\frac{1}{3} z^{3}_{011}\wedge z^{4}_{031} \wedge z^{5}_{113}
-\frac{1}{2} z^{3}_{011}\wedge z^{4}_{112}\wedge z^{5}_{121}
+z^{3}_{011}\wedge z^{4}_{102}\wedge z^{5}_{220} 
-z^{3}_{010}\wedge z^{4}_{210}\wedge z^{5}_{202}
-\frac{1}{2} z^{3}_{201}\wedge z^{4}_{102}\wedge z^{5}_{030}\\& 
-\frac{1}{2} z^{3}_{011}\wedge z^{4}_{002}\wedge z^{5}_{320} 
-\frac{1}{2} z^{3}_{011}\wedge z^{4}_{101}\wedge z^{5}_{221}
-\frac{1}{6} z^{3}_{011}\wedge z^{4}_{010}\wedge z^{5}_{401}
+\frac{1}{6} z^{3}_{011}\wedge z^{4}_{130}\wedge z^{5}_{014} 
-\frac{1}{2} z^{3}_{011}\wedge z^{4}_{202}\wedge z^{5}_{120}\\& 
+\frac{1}{2} z^{3}_{011}\wedge z^{4}_{103}\wedge z^{5}_{130}
+z^{3}_{011}\wedge z^{4}_{201}\wedge z^{5}_{121} 
-\frac{1}{4} z^{3}_{021}\wedge z^{4}_{211}\wedge z^{5}_{012}
-\frac{1}{2} z^{3}_{011}\wedge z^{4}_{030}\wedge z^{5}_{203}
+\frac{1}{2} z^{3}_{011}\wedge z^{4}_{021}\wedge z^{5}_{212} \\& 
-\frac{1}{2} z^{3}_{210}\wedge z^{4}_{201}\wedge z^{5}_{011}
+z^{3}_{011}\wedge z^{4}_{111}\wedge z^{5}_{122}
-\frac{1}{4} z^{3}_{010}\wedge z^{4}_{202}\wedge z^{5}_{210} 
+\frac{1}{2} z^{3}_{011}\wedge z^{4}_{220}\wedge z^{5}_{102}
-\frac{1}{2} z^{3}_{011}\wedge z^{4}_{003}\wedge z^{5}_{230}\\& 
+\frac{1}{4} z^{3}_{200}\wedge z^{4}_{301}\wedge z^{5}_{010} 
+\frac{1}{4} z^{3}_{010}\wedge z^{4}_{211}\wedge z^{5}_{201}  


\begin{align*} 
Q_1=& 1/2\Z{002}{3}\wedge\Z{003}{3}\wedge\Z{010}{3}\wedge\Z{030}{3}\wedge\Z{120}{3}\wedge\Z{201}{3}
+\Z{000}{3}\wedge\Z{011}{3}\wedge\Z{020}{3}\wedge\Z{102}{3}\wedge\Z{111}{3}\wedge\Z{300}{3} \\&
+3\Z{001}{3}\wedge\Z{011}{3}\wedge\Z{100}{3}\wedge\Z{111}{3}\wedge\Z{200}{3}\wedge\Z{210}{3}
-6\Z{001}{3}\wedge\Z{011}{3}\wedge\Z{101}{3}\wedge\Z{110}{3}\wedge\Z{111}{3}\wedge\Z{210}{3} \\&
-3/2\Z{001}{3}\wedge\Z{012}{3}\wedge\Z{020}{3}\wedge\Z{021}{3}\wedge\Z{102}{3}\wedge\Z{210}{3}
-3/2\Z{001}{3}\wedge\Z{012}{3}\wedge\Z{021}{3}\wedge\Z{100}{3}\wedge\Z{200}{3}\wedge\Z{210}{3} \\&
-3\Z{001}{3}\wedge\Z{012}{3}\wedge\Z{021}{3}\wedge\Z{101}{3}\wedge\Z{110}{3}\wedge\Z{210}{3}
+3/2\Z{001}{3}\wedge\Z{010}{3}\wedge\Z{020}{3}\wedge\Z{102}{3}\wedge\Z{201}{3}\wedge\Z{210}{3} \\&
-3/2\Z{001}{3}\wedge\Z{010}{3}\wedge\Z{100}{3}\wedge\Z{200}{3}\wedge\Z{201}{3}\wedge\Z{210}{3}
-3\Z{001}{3}\wedge\Z{010}{3}\wedge\Z{101}{3}\wedge\Z{110}{3}\wedge\Z{201}{3}\wedge\Z{210}{3} \\&
-3\Z{010}{3}\wedge\Z{011}{3}\wedge\Z{020}{3}\wedge\Z{102}{3}\wedge\Z{111}{3}\wedge\Z{201}{3}
+\Z{000}{3}\wedge\Z{002}{3}\wedge\Z{011}{3}\wedge\Z{111}{3}\wedge\Z{120}{3}\wedge\Z{300}{3} \\&
+1/2\Z{000}{3}\wedge\Z{002}{3}\wedge\Z{012}{3}\wedge\Z{021}{3}\wedge\Z{120}{3}\wedge\Z{300}{3}
-\Z{003}{3}\wedge\Z{011}{3}\wedge\Z{030}{3}\wedge\Z{100}{3}\wedge\Z{111}{3}\wedge\Z{200}{3} \\&
+3\Z{002}{3}\wedge\Z{100}{3}\wedge\Z{101}{3}\wedge\Z{110}{3}\wedge\Z{120}{3}\wedge\Z{200}{3}
-\Z{003}{3}\wedge\Z{010}{3}\wedge\Z{011}{3}\wedge\Z{030}{3}\wedge\Z{111}{3}\wedge\Z{201}{3} \\&
-3\Z{011}{3}\wedge\Z{020}{3}\wedge\Z{100}{3}\wedge\Z{102}{3}\wedge\Z{111}{3}\wedge\Z{200}{3}
+3\Z{011}{3}\wedge\Z{012}{3}\wedge\Z{020}{3}\wedge\Z{021}{3}\wedge\Z{102}{3}\wedge\Z{111}{3} \\&
-3/2\Z{002}{3}\wedge\Z{010}{3}\wedge\Z{020}{3}\wedge\Z{102}{3}\wedge\Z{120}{3}\wedge\Z{201}{3}
+2\Z{003}{3}\wedge\Z{011}{3}\wedge\Z{030}{3}\wedge\Z{101}{3}\wedge\Z{110}{3}\wedge\Z{111}{3} \\&
-1/2\Z{003}{3}\wedge\Z{012}{3}\wedge\Z{021}{3}\wedge\Z{030}{3}\wedge\Z{100}{3}\wedge\Z{200}{3}
-\Z{003}{3}\wedge\Z{012}{3}\wedge\Z{021}{3}\wedge\Z{030}{3}\wedge\Z{101}{3}\wedge\Z{110}{3} \\&
-3/2\Z{001}{3}\wedge\Z{020}{3}\wedge\Z{100}{3}\wedge\Z{102}{3}\wedge\Z{200}{3}\wedge\Z{210}{3}
-3\Z{001}{3}\wedge\Z{020}{3}\wedge\Z{101}{3}\wedge\Z{102}{3}\wedge\Z{110}{3}\wedge\Z{210}{3} \\&
-3\Z{002}{3}\wedge\Z{011}{3}\wedge\Z{100}{3}\wedge\Z{111}{3}\wedge\Z{120}{3}\wedge\Z{200}{3}
+\Z{003}{3}\wedge\Z{011}{3}\wedge\Z{020}{3}\wedge\Z{030}{3}\wedge\Z{102}{3}\wedge\Z{111}{3} \\&
-3/2\Z{002}{3}\wedge\Z{012}{3}\wedge\Z{020}{3}\wedge\Z{021}{3}\wedge\Z{102}{3}\wedge\Z{120}{3}
+3\Z{002}{3}\wedge\Z{011}{3}\wedge\Z{020}{3}\wedge\Z{102}{3}\wedge\Z{111}{3}\wedge\Z{120}{3} \\&
-\Z{002}{3}\wedge\Z{003}{3}\wedge\Z{011}{3}\wedge\Z{030}{3}\wedge\Z{111}{3}\wedge\Z{120}{3}
+6\Z{011}{3}\wedge\Z{012}{3}\wedge\Z{021}{3}\wedge\Z{101}{3}\wedge\Z{110}{3}\wedge\Z{111}{3} \\&
+3\Z{010}{3}\wedge\Z{100}{3}\wedge\Z{101}{3}\wedge\Z{110}{3}\wedge\Z{200}{3}\wedge\Z{201}{3}
-6\Z{011}{3}\wedge\Z{100}{3}\wedge\Z{101}{3}\wedge\Z{110}{3}\wedge\Z{111}{3}\wedge\Z{200}{3} \\&
-3/2\Z{001}{3}\wedge\Z{002}{3}\wedge\Z{010}{3}\wedge\Z{120}{3}\wedge\Z{201}{3}\wedge\Z{210}{3}
-\Z{000}{3}\wedge\Z{100}{3}\wedge\Z{101}{3}\wedge\Z{110}{3}\wedge\Z{200}{3}\wedge\Z{300}{3} \\&
-1/2\Z{001}{3}\wedge\Z{002}{3}\wedge\Z{003}{3}\wedge\Z{030}{3}\wedge\Z{120}{3}\wedge\Z{210}{3}
-3/2\Z{001}{3}\wedge\Z{002}{3}\wedge\Z{020}{3}\wedge\Z{102}{3}\wedge\Z{120}{3}\wedge\Z{210}{3} \\&
+3/2\Z{001}{3}\wedge\Z{002}{3}\wedge\Z{012}{3}\wedge\Z{021}{3}\wedge\Z{120}{3}\wedge\Z{210}{3}
+3\Z{001}{3}\wedge\Z{002}{3}\wedge\Z{011}{3}\wedge\Z{111}{3}\wedge\Z{120}{3}\wedge\Z{210}{3} \\&
-1/6\Z{000}{3}\wedge\Z{003}{3}\wedge\Z{010}{3}\wedge\Z{030}{3}\wedge\Z{201}{3}\wedge\Z{300}{3}
-1/2\Z{000}{3}\wedge\Z{002}{3}\wedge\Z{100}{3}\wedge\Z{120}{3}\wedge\Z{200}{3}\wedge\Z{300}{3} \\&
+\Z{000}{3}\wedge\Z{002}{3}\wedge\Z{101}{3}\wedge\Z{110}{3}\wedge\Z{120}{3}\wedge\Z{300}{3}
-1/2\Z{000}{3}\wedge\Z{002}{3}\wedge\Z{020}{3}\wedge\Z{102}{3}\wedge\Z{120}{3}\wedge\Z{300}{3} \\&
-\Z{003}{3}\wedge\Z{030}{3}\wedge\Z{100}{3}\wedge\Z{101}{3}\wedge\Z{110}{3}\wedge\Z{200}{3}
-3\Z{011}{3}\wedge\Z{012}{3}\wedge\Z{021}{3}\wedge\Z{100}{3}\wedge\Z{111}{3}\wedge\Z{200}{3} \\&
+1/2\Z{003}{3}\wedge\Z{010}{3}\wedge\Z{012}{3}\wedge\Z{021}{3}\wedge\Z{030}{3}\wedge\Z{201}{3}
-3\Z{020}{3}\wedge\Z{100}{3}\wedge\Z{101}{3}\wedge\Z{102}{3}\wedge\Z{110}{3}\wedge\Z{200}{3} \\&
+1/2\Z{002}{3}\wedge\Z{003}{3}\wedge\Z{012}{3}\wedge\Z{021}{3}\wedge\Z{030}{3}\wedge\Z{120}{3}
+3\Z{001}{3}\wedge\Z{011}{3}\wedge\Z{020}{3}\wedge\Z{102}{3}\wedge\Z{111}{3}\wedge\Z{210}{3} \\&
-3/2\Z{012}{3}\wedge\Z{020}{3}\wedge\Z{021}{3}\wedge\Z{100}{3}\wedge\Z{102}{3}\wedge\Z{200}{3}
+1/2\Z{000}{3}\wedge\Z{001}{3}\wedge\Z{012}{3}\wedge\Z{021}{3}\wedge\Z{210}{3}\wedge\Z{300}{3} \\&
-1/2\Z{000}{3}\wedge\Z{001}{3}\wedge\Z{020}{3}\wedge\Z{102}{3}\wedge\Z{210}{3}\wedge\Z{300}{3}
-1/2\Z{000}{3}\wedge\Z{001}{3}\wedge\Z{002}{3}\wedge\Z{120}{3}\wedge\Z{210}{3}\wedge\Z{300}{3} \\&
-1/6\Z{000}{3}\wedge\Z{001}{3}\wedge\Z{003}{3}\wedge\Z{030}{3}\wedge\Z{210}{3}\wedge\Z{300}{3}
+1/2\Z{000}{3}\wedge\Z{001}{3}\wedge\Z{010}{3}\wedge\Z{201}{3}\wedge\Z{210}{3}\wedge\Z{300}{3} \\&
+\Z{000}{3}\wedge\Z{001}{3}\wedge\Z{011}{3}\wedge\Z{111}{3}\wedge\Z{210}{3}\wedge\Z{300}{3}
-1/3\Z{000}{3}\wedge\Z{003}{3}\wedge\Z{011}{3}\wedge\Z{030}{3}\wedge\Z{111}{3}\wedge\Z{300}{3} \\&
-1/2\Z{000}{3}\wedge\Z{020}{3}\wedge\Z{100}{3}\wedge\Z{102}{3}\wedge\Z{200}{3}\wedge\Z{300}{3}
-\Z{000}{3}\wedge\Z{020}{3}\wedge\Z{101}{3}\wedge\Z{102}{3}\wedge\Z{110}{3}\wedge\Z{300}{3} \\&
+\Z{003}{3}\wedge\Z{010}{3}\wedge\Z{030}{3}\wedge\Z{101}{3}\wedge\Z{110}{3}\wedge\Z{201}{3}
+6\Z{010}{3}\wedge\Z{011}{3}\wedge\Z{101}{3}\wedge\Z{110}{3}\wedge\Z{111}{3}\wedge\Z{201}{3} \\&
+3/2\Z{010}{3}\wedge\Z{020}{3}\wedge\Z{100}{3}\wedge\Z{102}{3}\wedge\Z{200}{3}\wedge\Z{201}{3}
+3\Z{010}{3}\wedge\Z{020}{3}\wedge\Z{101}{3}\wedge\Z{102}{3}\wedge\Z{110}{3}\wedge\Z{201}{3} \\&
+3/2\Z{002}{3}\wedge\Z{010}{3}\wedge\Z{012}{3}\wedge\Z{021}{3}\wedge\Z{120}{3}\wedge\Z{201}{3}
+3\Z{002}{3}\wedge\Z{010}{3}\wedge\Z{011}{3}\wedge\Z{111}{3}\wedge\Z{120}{3}\wedge\Z{201}{3} \\&
-1/2\Z{002}{3}\wedge\Z{003}{3}\wedge\Z{030}{3}\wedge\Z{100}{3}\wedge\Z{120}{3}\wedge\Z{200}{3}
+\Z{002}{3}\wedge\Z{003}{3}\wedge\Z{030}{3}\wedge\Z{101}{3}\wedge\Z{110}{3}\wedge\Z{120}{3} \\&
+1/2\Z{003}{3}\wedge\Z{012}{3}\wedge\Z{020}{3}\wedge\Z{021}{3}\wedge\Z{030}{3}\wedge\Z{102}{3}
+1/2\Z{002}{3}\wedge\Z{003}{3}\wedge\Z{020}{3}\wedge\Z{030}{3}\wedge\Z{102}{3}\wedge\Z{120}{3} \\&
-1/6\Z{000}{3}\wedge\Z{002}{3}\wedge\Z{003}{3}\wedge\Z{030}{3}\wedge\Z{120}{3}\wedge\Z{300}{3}
-1/2\Z{000}{3}\wedge\Z{002}{3}\wedge\Z{010}{3}\wedge\Z{120}{3}\wedge\Z{201}{3}\wedge\Z{300}{3} \\&
+\Z{000}{3}\wedge\Z{001}{3}\wedge\Z{101}{3}\wedge\Z{110}{3}\wedge\Z{210}{3}\wedge\Z{300}{3}
+3\Z{012}{3}\wedge\Z{021}{3}\wedge\Z{100}{3}\wedge\Z{101}{3}\wedge\Z{110}{3}\wedge\Z{200}{3} \\&
+3\Z{010}{3}\wedge\Z{011}{3}\wedge\Z{012}{3}\wedge\Z{021}{3}\wedge\Z{111}{3}\wedge\Z{201}{3}
-3/2\Z{002}{3}\wedge\Z{010}{3}\wedge\Z{100}{3}\wedge\Z{120}{3}\wedge\Z{200}{3}\wedge\Z{201}{3} \\&
+3\Z{002}{3}\wedge\Z{010}{3}\wedge\Z{101}{3}\wedge\Z{110}{3}\wedge\Z{120}{3}\wedge\Z{201}{3}
+1/2\Z{003}{3}\wedge\Z{010}{3}\wedge\Z{020}{3}\wedge\Z{030}{3}\wedge\Z{102}{3}\wedge\Z{201}{3} \\&
-3\Z{002}{3}\wedge\Z{011}{3}\wedge\Z{012}{3}\wedge\Z{021}{3}\wedge\Z{111}{3}\wedge\Z{120}{3}
+1/6\Z{000}{3}\wedge\Z{003}{3}\wedge\Z{020}{3}\wedge\Z{030}{3}\wedge\Z{102}{3}\wedge\Z{300}{3} \\&
+1/6\Z{000}{3}\wedge\Z{003}{3}\wedge\Z{012}{3}\wedge\Z{021}{3}\wedge\Z{030}{3}\wedge\Z{300}{3}
+1/2\Z{003}{3}\wedge\Z{020}{3}\wedge\Z{030}{3}\wedge\Z{100}{3}\wedge\Z{102}{3}\wedge\Z{200}{3} \\&
-1/2\Z{000}{3}\wedge\Z{012}{3}\wedge\Z{021}{3}\wedge\Z{100}{3}\wedge\Z{200}{3}\wedge\Z{300}{3}
-\Z{000}{3}\wedge\Z{012}{3}\wedge\Z{021}{3}\wedge\Z{101}{3}\wedge\Z{110}{3}\wedge\Z{300}{3} \\&
-1/2\Z{000}{3}\wedge\Z{012}{3}\wedge\Z{020}{3}\wedge\Z{021}{3}\wedge\Z{102}{3}\wedge\Z{300}{3}
+\Z{003}{3}\wedge\Z{020}{3}\wedge\Z{030}{3}\wedge\Z{101}{3}\wedge\Z{102}{3}\wedge\Z{110}{3} \\&
-2\Z{000}{3}\wedge\Z{011}{3}\wedge\Z{101}{3}\wedge\Z{110}{3}\wedge\Z{111}{3}\wedge\Z{300}{3}
-\Z{000}{3}\wedge\Z{011}{3}\wedge\Z{012}{3}\wedge\Z{021}{3}\wedge\Z{111}{3}\wedge\Z{300}{3} \\&
-3\Z{010}{3}\wedge\Z{011}{3}\wedge\Z{100}{3}\wedge\Z{111}{3}\wedge\Z{200}{3}\wedge\Z{201}{3}
-3/2\Z{001}{3}\wedge\Z{002}{3}\wedge\Z{100}{3}\wedge\Z{120}{3}\wedge\Z{200}{3}\wedge\Z{210}{3} \\&
-1/2\Z{000}{3}\wedge\Z{010}{3}\wedge\Z{100}{3}\wedge\Z{200}{3}\wedge\Z{201}{3}\wedge\Z{300}{3}
-\Z{000}{3}\wedge\Z{010}{3}\wedge\Z{101}{3}\wedge\Z{110}{3}\wedge\Z{201}{3}\wedge\Z{300}{3} \\&
-1/2\Z{000}{3}\wedge\Z{010}{3}\wedge\Z{012}{3}\wedge\Z{021}{3}\wedge\Z{201}{3}\wedge\Z{300}{3}
+1/2\Z{000}{3}\wedge\Z{010}{3}\wedge\Z{020}{3}\wedge\Z{102}{3}\wedge\Z{201}{3}\wedge\Z{300}{3} \\&
-\Z{000}{3}\wedge\Z{010}{3}\wedge\Z{011}{3}\wedge\Z{111}{3}\wedge\Z{201}{3}\wedge\Z{300}{3}
+1/6\Z{000}{3}\wedge\Z{003}{3}\wedge\Z{030}{3}\wedge\Z{100}{3}\wedge\Z{200}{3}\wedge\Z{300}{3} \\&
+1/3\Z{000}{3}\wedge\Z{003}{3}\wedge\Z{030}{3}\wedge\Z{101}{3}\wedge\Z{110}{3}\wedge\Z{300}{3}
-3\Z{001}{3}\wedge\Z{010}{3}\wedge\Z{011}{3}\wedge\Z{111}{3}\wedge\Z{201}{3}\wedge\Z{210}{3} \\&
-3/2\Z{001}{3}\wedge\Z{010}{3}\wedge\Z{012}{3}\wedge\Z{021}{3}\wedge\Z{201}{3}\wedge\Z{210}{3}
+1/2\Z{001}{3}\wedge\Z{003}{3}\wedge\Z{030}{3}\wedge\Z{100}{3}\wedge\Z{200}{3}\wedge\Z{210}{3} \\&
+\Z{001}{3}\wedge\Z{003}{3}\wedge\Z{030}{3}\wedge\Z{101}{3}\wedge\Z{110}{3}\wedge\Z{210}{3}
+3/2\Z{010}{3}\wedge\Z{012}{3}\wedge\Z{021}{3}\wedge\Z{100}{3}\wedge\Z{200}{3}\wedge\Z{201}{3} \\&
+1/2\Z{001}{3}\wedge\Z{003}{3}\wedge\Z{020}{3}\wedge\Z{030}{3}\wedge\Z{102}{3}\wedge\Z{210}{3}
+\Z{003}{3}\wedge\Z{011}{3}\wedge\Z{012}{3}\wedge\Z{021}{3}\wedge\Z{030}{3}\wedge\Z{111}{3} \\&
-3\Z{001}{3}\wedge\Z{011}{3}\wedge\Z{012}{3}\wedge\Z{021}{3}\wedge\Z{111}{3}\wedge\Z{210}{3}
+1/2\Z{003}{3}\wedge\Z{010}{3}\wedge\Z{030}{3}\wedge\Z{100}{3}\wedge\Z{200}{3}\wedge\Z{201}{3} \\&
+1/2\Z{001}{3}\wedge\Z{003}{3}\wedge\Z{012}{3}\wedge\Z{021}{3}\wedge\Z{030}{3}\wedge\Z{210}{3}
+3\Z{010}{3}\wedge\Z{012}{3}\wedge\Z{021}{3}\wedge\Z{101}{3}\wedge\Z{110}{3}\wedge\Z{201}{3} \\&
-\Z{001}{3}\wedge\Z{003}{3}\wedge\Z{011}{3}\wedge\Z{030}{3}\wedge\Z{111}{3}\wedge\Z{210}{3}
+3/2\Z{010}{3}\wedge\Z{012}{3}\wedge\Z{020}{3}\wedge\Z{021}{3}\wedge\Z{102}{3}\wedge\Z{201}{3} \\&
+3\Z{001}{3}\wedge\Z{002}{3}\wedge\Z{101}{3}\wedge\Z{110}{3}\wedge\Z{120}{3}\wedge\Z{210}{3}
-1/2\Z{001}{3}\wedge\Z{003}{3}\wedge\Z{010}{3}\wedge\Z{030}{3}\wedge\Z{201}{3}\wedge\Z{210}{3} \\&
+\Z{000}{3}\wedge\Z{011}{3}\wedge\Z{100}{3}\wedge\Z{111}{3}\wedge\Z{200}{3}\wedge\Z{300}{3}
+3/2\Z{002}{3}\wedge\Z{012}{3}\wedge\Z{021}{3}\wedge\Z{100}{3}\wedge\Z{120}{3}\wedge\Z{200}{3} \\&
-3\Z{002}{3}\wedge\Z{012}{3}\wedge\Z{021}{3}\wedge\Z{101}{3}\wedge\Z{110}{3}\wedge\Z{120}{3}
+3/2\Z{002}{3}\wedge\Z{020}{3}\wedge\Z{100}{3}\wedge\Z{102}{3}\wedge\Z{120}{3}\wedge\Z{200}{3} \\&
-3\Z{002}{3}\wedge\Z{020}{3}\wedge\Z{101}{3}\wedge\Z{102}{3}\wedge\Z{110}{3}\wedge\Z{120}{3}
-3\Z{012}{3}\wedge\Z{020}{3}\wedge\Z{021}{3}\wedge\Z{101}{3}\wedge\Z{102}{3}\wedge\Z{110}{3} \\&
-6\Z{002}{3}\wedge\Z{011}{3}\wedge\Z{101}{3}\wedge\Z{110}{3}\wedge\Z{111}{3}\wedge\Z{120}{3}
-3\Z{001}{3}\wedge\Z{100}{3}\wedge\Z{101}{3}\wedge\Z{110}{3}\wedge\Z{200}{3}\wedge\Z{210}{3} 
\end{align*}



}%
\end{document}